\documentclass[1p,authoryear]{elsarticle}

\usepackage{amsfonts}
\usepackage{amsmath}
\usepackage{amssymb}
\usepackage[applemac]{inputenc}
\usepackage{calrsfs} 

\usepackage{graphicx}
\usepackage{subfigure}
\usepackage{mathabx}

\usepackage{natbib}
\usepackage{geometry}
\usepackage{color}

\newcommand{\vs}[1]{\vspace{#1}}    
\newcommand{\hs}[1]{\hspace{#1}}    
\newcommand{\vset}[1]{\vspace*{#1}}    
    
\newcommand{\begitem}{\begin{itemize}}    
\newcommand{\finit}{\end{itemize}}    
\newcommand{\begenum}{\begin{enumerate}}    
\newcommand{\finenum}{\end{enumerate}}    
\newcommand{\begar}{\begin{array}}    
\newcommand{\finar}{\end{array}}    
\newcommand{\begeq}[1]{\begin{equation} \label{#1}}    
\newcommand{\fineq}{\end{equation}}    
\newcommand{\begct}{\begin{center}}    
\newcommand{\finct}{\end{center}}

\newcommand{\zun}{\vs{0.1cm}} 
     
\newcommand{\zdeux}{\vs{0.2cm}} 
\newcommand{\ztrois}{\vs{0.3cm}}    
\newcommand{\zcinq}{\vs{0.5cm}} 
    
\newcommand{\tinf}{\rightarrow\infty}         
\newcommand{\tqdninf}{\stackrel{n\tinf}{\longrightarrow}} 

\newcommand{\bE}{\mathbb{E}}

\newcommand{\bV}{\mathbb{V}}

\newcommand{\bP}{\mathbb{P}}

\newcommand{\bI}{\mathbb{I}}

\newcommand{\limloi}{ \stackrel{d}{\longrightarrow} }
\newcommand{\limprob}{ \stackrel{\bP}{\longrightarrow} }
\newcommand{\equloi}{ \stackrel{d}{=} }

\newcommand{\Fbar}{\widebar{F}}
\newcommand{\Fbarn}{\widebar{F}^{KM}_n}

\newcommand{\Hbar}{\widebar{H}}
 
\newcommand{\Znmiun}{{Z_{n-i+1,n}}}

\newcommand{\Znmjun}{{Z_{n-j+1,n}}}
\newcommand{\Znmjunmin}{{Z^-_{n-j+1,n}}}
\newcommand{\Znmj}{{Z_{n-j,n}}}
\newcommand{\Znmk}{{Z_{n-k,n}}}
\newcommand{\Ytildekmjun}{{\tilde Y_{k-j+1,k}}}
\newcommand{\Ytildekmjunexpo}[1]{{\tilde Y_{k-j+1,k}^{#1}}}

\newcommand{\Znmjunbetaun}{{Z_{n-j+1,n}^{-\beta_1}}}
\newcommand{\Znmkbetaun}{{Z_{n-k,n}^{-\beta_1}}}
 
\newcommand{\RFj}{{RF_{j}}}  
\newcommand{\RFjbet}{{RF_{j,\beta}}}
\newcommand{\RFchapj}{{\widehat{RF}_{j}}}  
\newcommand{\RapFj}{{ \frac{ \Fbar(\Znmjun) }{ \Fbar(\Znmk) } }}  
\newcommand{\RapFchapj}{{ \frac{ \Fbarn(\Znmjun) }{ \Fbarn(\Znmk) } }}  
\newcommand{\RapFchapjmin}{{ \frac{ \Fbarn(\Znmjunmin) }{ \Fbarn(\Znmk) } }} 

\newcommand{\ksi}{\xi}
\newcommand{\ksijkbet}{\xi_{j,k,\beta}}
\newcommand{\ksijkbetp}{\xi'_{j,k,\beta}}

\newcommand{\ksipjsurj}{\frac{\xi'_j}{j}}
\newcommand{\ksijdev}{j\log\frac{\Znmjun}{\Znmj}}

\newcommand{\ui}{u_{i,k}}
\newcommand{\uiplusun}{u_{i+1,k}}

\newcommand{\uipmunbet}{u_{i,k}^{\pbet-1}}

\newcommand{\uj}{u_{j,k}}
\newcommand{\ujmoinsun}{u_{j-1,k}}

\newcommand{\ujpmunbet}{u_{j,k}^{\pbet-1}}

\newcommand{\ujpbet}{u_{j,k,}^{\pbet}}

\newcommand{\Vjk}{V_{j,k}}

\newcommand{\Vjkpbet}{V_{j,k}^{\pbet}}

\newcommand{\calEn}{{{\mathcal E}_{n,\alpha}}}

\newcommand{\deltanmiun}{\delta_{n-i+1,n}}

\newcommand{\bnk}{b_{n,k}}
\newcommand{\Rnj}{R_{n,j}}
\newcommand{\Rnjbet}{R_{n,j,\beta}}
\newcommand{\Tkn}[1]{T_{k,n}^{(#1)}}
\newcommand{\Sik}{S_{i,k}}
\newcommand{\Sikbet}{S_{i,k,\beta}}

\newcommand{\Ein}{E^{(n)}_i}
\newcommand{\Eln}{E^{(n)}_l}
\newcommand{\Ejn}{E^{(n)}_j}

\newcommand{\sumideuxk}{{\sum_{i=2}^k}}
\newcommand{\sumjdeuxk}{{\sum_{j=2}^k}}
\newcommand{\sumideuxj}{{\sum_{i=2}^j}}
\newcommand{\sumjdeuxi}{{\sum_{j=2}^i}}

\newcommand{\summ}[2]{{ \sum_{#1}^{#2} }}

\newcommand{\gamchapunk}{\widehat\gamma^{(W)}_{1,k}}
\newcommand{\gamchabeta}{\widehat\gamma_{1,k}(\beta)}
\newcommand{\gamBR}{\widehat\gamma^{(BR)}_{1,k}}
\newcommand{\TB}{\widehat{T}_{k}(\beta)}

\newcommand{\gamchapHillunk}{\widehat\gamma^{(H)}_{1,k}}

\newcommand{\bet}{b_*}
\newcommand{\btildeet}{\tilde b_*}

\newcommand{\pbet}{p_{\beta}}
\newcommand{\mT}{m_{\beta}}
\newcommand{\sigT}{\sigma_{\beta}}
\newcommand{\mg}{m_{\gamma_1,\beta}}
\newcommand{\sigg}{\sigma_{\gamma_1,\beta}}

\newcommand{\sigBR}{\sigma_{(BR)}}

\newcommand{\oPdeun}{o_{\bP}(1)}

\newcommand{\gbetaet}{\gamma\beta_*}

\begin{document}

\newtheorem{theo}{Theorem}
\newtheorem{prop}{Proposition}
\newtheorem{defi}{Definition}
\newtheorem{lem}{Lemma}
\newtheorem{cor}{Corollary}
\newtheorem{rmk}{Remark}
  
\begin{center}
{\Large {\sc Estimation of the extreme value index  in a censorship framework: asymptotic and finite sample behaviour}} 
\bigskip

\large Jan Beirlant(1) ,  Julien Worms\footnote{Corresponding author J. Worms} (2) \& Rym Worms (3)

%

\vset{1.cm}
 (1) KU Leuven \\ 
  Department of Mathematics and Leuven Statistics Research Center, KU Leuven,  Belgium, \\ Department of Mathematical Statistics and Actuarial Science, University of the Free State, South Africa \\
   e-mail : {\tt jan.beirlant@kuleuven.be}
\bigskip

\vset{1.cm}
 (2) Universit\'e de Versailles-Saint-Quentin-En-Yvelines\\
 Laboratoire de Math\'ematiques de Versailles (CNRS UMR 8100), \\
 F-78035 Versailles Cedex, France, \\
 e-mail : {\tt julien.worms@uvsq.fr}
\bigskip
\vset{1.cm}

(3) Universit\'e Paris-Est \\
Laboratoire d'Analyse et de Math\'ematiques Appliqu\'ees \\
(CNRS UMR 8050), \\
 UPEMLV, UPEC, F-94010, Cr\'eteil, France, \\
 e-mail : {\tt rym.worms@u-pec.fr}
\end{center}
\vspace{1.cm}

\begin{center}
{\bf Abstract}
\end{center}
 We revisit the estimation of the extreme value index for randomly censored data from a heavy tailed distribution. We introduce a new class of estimators which encompasses earlier proposals given in Worms and Worms (2014) and Beirlant et al. (2018), which were shown to have good bias properties compared with the pseudo maximum likelihood estimator proposed in Beirlant et al. (2007) and Einmahl et al. (2008). However the asymptotic normality of the type of estimators first proposed in Worms and Worms (2014) was still lacking. We  derive an asymptotic representation and the asymptotic normality of the larger class of estimators  and consider their  finite sample behaviour.  Special attention is paid to the case of heavy censoring, i.e. where the amount of censoring in the tail is at least 50\%. We obtain the asymptotic normality with a classical $\sqrt{k}$ rate where $k$ denotes the number of top data used in the estimation, depending on the degree of censoring.   

\vfill

\noindent 
{\it AMS Classification. } Primary 62G32 ; Secondary 62N02 
  \vspace{0.1cm} \\
{\it Keywords and phrases.~} Extreme value index. Tail inference. Random censoring.  Asymptotic representation.  \\

\setlength{\belowdisplayskip}{0.2cm}
\setlength{\abovedisplayskip}{0.2cm}

\newpage

\section{Introduction}
 \label{intro}

Starting from  \citet{Beirlant07}, the estimation of the extreme value index  in a censorship framework  is of growing interest.   Suppose we observe a sample of $n$ independent couples $(Z_i,\delta_i)_{1\leq i\leq n}$ where
\[
 Z_i = \min(X_i,C_i) \hs{0.3cm}  \mbox{and} \hs{0.3cm} \delta_i=\bI_{X_i\leq C_i}. 
 \]
The i.i.d. samples $(X_i)_{i\leq n}$ and $(C_i)_{i\leq n}$, of respective continuous distribution functions $F$ and $G$, are samples from the variable of interest $X$ and of the censoring variable $C$, measured on $n$ individual items (insurance claims, hospitalized patients, ...). The variables $X$ and $C$ are supposed to be independent and, for convenience only, we will suppose in this work that they are non-negative. We will denote by $Z_{1,n}\leq \ldots \leq Z_{i,n}\leq \ldots \leq Z_{n,n}$ the order statistics associated to the observed sample, and by $(\delta_{1,n},\ldots,\delta_{n,n})$ the corresponding indicators of non-censorship.
\zdeux 

\citet{Einmahl08} presented a general method for adapting estimators of the extreme value index in this censorship framework.  \citet{WormsWorms14} proposed a more survival analysis-oriented approach restricted to the heavy tail case, while \citet{Diop2014} extended the framework to data with covariate information. \citet{BeirlantBardoutsosWetGijbels16} and \citet{BeirlantMaribeVester18} proposed bias-reduced versions of two existing estimators. See also \citet{Necir2015}, \citet{Necir2016} and \citet{Necir2018} for other papers on the subject. 

\zdeux 

In this paper, we propose a new class of estimators that encompasses one of the estimators proposed in \citet{WormsWorms14}  and propose a novel approach to prove the asymptotic normality of these estimators which was unknown up to now for the case $\beta=0$.  We consider here that the distributions $F$ and $G$ are heavy-tailed, with positive and respective extreme value indices (EVI) $\gamma_1$ and $\gamma_2$, i.e.
\[
\bar{F}(x)=  1-F(x)  = x^{-1/\gamma_1} l_F(x)  \mbox{ \ and \ }  \bar{G}(y)=  1-G(y)=  y^{-1/\gamma_2} l_G(y), 
\]
where $l_F$ and $l_G$ are slowly varying at infinity. Our target is the EVI $\gamma_1$, which we try to recover from our randomly censored observations.
\zun

Denoting the distribution function of $Z$ with $H$, by independence of $X$ and $C$ we readily obtain $\bar{H}(z)= 1-H(z)=  z^{-1/\gamma} l_H(z) $, 
where $l_H=l_F l_G$ and the EVI $\gamma$ of $Z$ is related to those of $X$ and $C$ via the important relation $1/\gamma= 1/\gamma_1 +1/\gamma_2$. Further in this paper, we will denote by $p$ the crucial quantity $p=\gamma/\gamma_1=\gamma_2/(\gamma_1+\gamma_2)\in ]0,1[$, which has to be interpreted as the asymptotic proportion of non-censored observations in the tail.
\zun 

 We assume in this work that  the slowly varying functions $l_F$ and $l_G$ satisfy the second order condition  first proposed  by \citet{HallWelsh85}. This yields the so called "Hall-type" model, i.e. as $x, y \rightarrow + \infty$, 
 \begin{align}
 \bar{F}(x) & = C_1  x^{-1/\gamma_1} \left(  1 + D_1 x^{-\beta_1} (1+o(1))  \right)   \label{condFbar} \\
 \bar{G}(y) & =  C_2  y^{-1/\gamma_2} \left(  1 + D_2 y^{-\beta_2} (1+o(1))  \right)   \label{condGbar}
\end{align}
where $\beta_1$, $\beta_2$, $C_1$,$C_2$ are positive constants and  $D_1$, $D_2$  are real constants. Then, setting  
\[
  C=C_1C_2, \hspace{0.4cm} \beta_*=\min(\beta_1,\beta_2),  \makebox[1.3cm][c]{ and }
  D_*= \left\{  \begar{lll}  D_1  &  \mbox{ if }  & \beta_1 < \beta_2, \\
                                                D_2   & \mbox{ if }  & \beta_2 < \beta_1, \\
						D_1 + D_2 & \mbox{ if }  & \beta_1 = \beta_2, 
						 \finar  \right.
\]
we have, as $z\to\infty$,
\begeq{condHbar}
  \bar{H}(z)=  C  z^{-1/\gamma} \left(  1 + D_* z^{-\beta_*} (1+o(1))  \right).
\fineq 
Correspondingly, with $H^-(u) =  \inf \{ z: H(z) \geq u \} $ ($0<u<1$) the quantile function corresponding to $H$, we consider $U_H(x)=H^-(1-1/x)$, the right-tail function of $H$, for which as $x\to \infty$, 
\begeq{condUH}
U_H(x)= C^{\gamma}  x^{\gamma} \left(  1 + \gamma D_* C^{-\beta_* \gamma }  x^{-\beta_* \gamma } (1+o(1))  \right).
\medskip
\fineq

Let us now explain how we build our new family of estimators of $\gamma_1$. For some real number $\beta$, consider the Box-Cox transform $k_{-\beta}(u) = \int_1^u t^{-\beta -1}dt$ for $u>1$, with the case $\beta=0$ leading to $k_0(u)=\log(u)$.  Based on the relation 
\begeq{basiclimit}
 \lim_{t\rightarrow\infty} \bE\, [ \,k_{-\beta} (X/t) \, | \, X>t \, ] =  
 \lim_{t\rightarrow\infty} \int_1^{\infty} \frac{\Fbar(ut)}{\Fbar(t)} dk_{-\beta}(u) = {\gamma_1 \over 1+\beta \gamma_1}, 
\fineq
and estimating $\Fbar$ by the Kaplan-Meier estimator $\Fbarn$   defined for $t<Z_{n,n}$ by
\begeq{KM}
 \Fbarn(t) = \prod_{Z_{i,n}\leq t} \left( \frac{n-i}{n-i+1} \right)^{\delta_{i,n}},
\fineq
we introduce the following class of statistics   
\begeq{defTB}
 \TB := \summ{j=2}{k} \RapFchapj \left( k_{-\beta}\left(\frac{\Znmjun}{\Znmk}\right) -  k_{-\beta}\left(\frac{\Znmj}{\Znmk}\right)  \right) 
\fineq
where $k=k_n$ denotes an integer sequence satisfying $k_n\tinf$ and $k_n=o(n)$. 
With $\beta=0$ we thus obtain the estimator 
\begeq{defgamchap1}
 \gamchapunk := \widehat{T}_k (0)=\summ{j=2}{k} \RapFchapj \log\left(\frac{\Znmjun}{\Znmj}\right)
\fineq
of $\gamma_1$ which was considered in Worms and Worms (2014) and Beirlant et al. (2018). In fact
$\gamchapunk$  turns out to be very close to the estimator
$ \summ{j=1}{k} \RapFchapjmin \log\frac{\Znmjun}{\Znmj}
$
 defined in equation (12) of \citet{WormsWorms14} based on ideas issued from the so-called Leurgans approach in survival regression analysis. The difference concerns a different way to circumvent the use of $\Fbarn$ at $Z_{n,n}$: whether using left-limits or deleting    $\Fbarn (Z_{n,n})$ as in $\gamchapunk$. 
 \medskip

Note that  the statistics $\TB$ were used in  \citet{BeirlantMaribeVester18} to obtain a bias-reduced version of the  estimator $\gamchapunk$ : 
\begeq{defEstimBRBeirlant}
\gamBR = \gamchapunk   -  \frac{(1+\beta_1\gamchapunk)^2 (1+2\beta_1\gamchapunk)}{(\beta_1\gamchapunk)^2}
\left(\widehat{T}_{k}(\beta_1) -\frac{\gamchapunk}{1+\beta_1\gamchapunk}\right),
\fineq
where $\beta_1$ denotes the second order parameter of $F$ in assumption  \eqref{condFbar}.
 \medskip
 
Now, it is clear from \eqref{basiclimit}  
that we can  construct the following  estimator of $\gamma_1$ when the tuning parameter $\beta$ is supposed to be larger than $-1/\gamma_1$: 
\begeq{defgambeta}
\gamchabeta = \frac{\TB}{1-\beta \TB}.
\fineq
\zdeux
We will compare these estimators with the pseudo maximum likelihood estimator which was first proposed in the random censoring context by \citet{Beirlant07} and \citet{Einmahl08}:   
\begeq{defEstimHillBeirlant}
 \gamchapHillunk = \frac 1{\widehat p_n} \, \frac 1 k \summ{i=1}{k}  \log\frac{\Znmiun}{\Znmk} 
 \makebox[1.5cm][c]{where} 
 \widehat p_n = \frac 1 k \summ{i=1}{k} \deltanmiun.
\smallskip
\fineq

In \citet{BeirlantMaribeVester18} a small sample simulation study was performed using all those available estimators and it was found that $\gamchapunk$ overall shows quite good bias and MSE performance. However, since no results on the asymptotic normality of this estimator were available yet, these authors proposed the use of a bootstrap algorithm to construct confidence intervals. In this paper we prove the  asymptotic normality of $\gamchabeta$ in the case $p+\beta\gamma>{1 \over 2}$. Hence this paper provides the first complete proof of the asymptotic normality  for $\gamchapunk$ in case $p>{1 \over 2}$, issued from an explicit asymptotic development stated in Theorem \ref{main-theorem} of the next section. In the  deterministic threshold case, this central limit result (for $\gamchapunk$) had already been obtained in \citet{WormsWorms18}, where a more general competing risks setting was considered, and using a different approach from the present proof. 
\smallskip

The restriction  $p>{1 \over 2}$ is rather restrictive for instance in insurance problems such as those discussed in \citet{BeirlantMaribeVester18} where heavy censoring appears. The introduction of the class of estimators $\gamchabeta$ helps to circumvent this problem when considering $\beta >0$. 
\smallskip

Finally, in the next section, we will see that our results also lead to the statement of the asymptotic normality of the bias-reduced estimator $\gamBR$, which was not known so far. 
\zdeux

\zdeux
Our paper is organized as follows: in Section \ref{Results}, we state and discuss the asymptotic normality result for $\gamchabeta$ and $\gamBR$. Section \ref{Proofs}  is devoted to the proof.  Technical aspects of the proof are postponed to the Appendix. In Section \ref{FiniteSample} we discuss the finite sample behavior of the different estimators  $\gamchabeta$ with $\beta >-1/\gamma_1$, and of $\gamBR$.
\zcinq

\section{Results}  \label{Results}

Our first and main result states the asymptotic behavior  of the statistics $\TB$ defined in (\ref{defTB}).  This result entails the asymptotic normality of the estimator $\gamchapunk$ of $\gamma_1$ by considering the particular case  $\beta=0$. The main condition is that  the heaviness of the tail of the censoring variable $C$ should be sufficiently  high with respect to the one of the variable $X$.
More precisely, introducing the notation $\pbet = p+\gamma \beta=p(1+\gamma_1\beta)$, the condition is be that $\pbet$ must be larger than $1/2$ ({\it i.e.} $\gamma_2>\gamma_1/(1+2\gamma_1\beta)$).

 \begin{theo} \label{main-theorem}
 Let conditions $(\ref{condFbar})$ and $(\ref{condGbar})$ hold. We assume further that $p_{\beta} > \frac 1 2$, and 
 \begeq{conditionbiais}
  \sqrt{k} \left( k/n\right)^{\gamma \beta_*} \; \tqdninf \; \lambda,
\fineq
and, if $\lambda=0$,   that $n=O(k_n^B)$ for some large enough $B>0$. We then have, as $n\tinf$, 
\[
 \sqrt{k} \left(\TB - \frac{\gamma_1}{1+\gamma_1 \beta} \right)  = G_n + \lambda \mT + o_{\bP}(1) \makebox[1.5cm][c]{where}  
 G_n \equloi \frac{\gamma}{\pbet} \frac{1}{\sqrt{k}} \sum_{i=2}^k  \ui^{\pbet -1} \left( p(E_i-1)- (\bI_{U_i \leq p}-p) \right) 
\]
with $(E_i)$ and $(U_i)$ denoting independent iid samples with, respectively, standard exponential and standard uniform distributions, and 
\[  
 \mT = \left\{ \begar[c]{ll} 
  -\gamma^2 \beta_1 D_1 C^{-\gamma \beta_1} \pbet^{-1} (\pbet+\gamma \beta_1)^{-1}  & \mbox{  if $\beta_1 \leq \beta_2 $},  \\
  0  & \mbox{   if $\beta_1 > \beta_2$.}
  \finar \right.
\]
Therefore, as $n\tinf$,
\[
 \sqrt{k} \left(\TB - \frac{\gamma_1}{1+\gamma_1 \beta} \right)  \limloi N(\lambda \mT,\sigT^2)
 \makebox[1.5cm][c]{where}
 \sigT^2=\displaystyle \frac{\gamma^2}{\pbet^2} \frac{p}{2 \pbet -1}
 =\displaystyle \gamma_1^2\frac{p}{2p-1} \frac{p^2}{\pbet^2} \frac{2p-1}{2 \pbet -1}.
\]
 \end{theo}

Since $G_n$ is a sum of independent random variables, it is then easy, using Lyapunov's CLT and the delta-method, to derive the following asymptotic normality result for the family of estimators $\gamchabeta$ of $\gamma_1$ defined by (\ref{defgambeta}). 

\begin{cor} \label{corGammaW}
Under the conditions of Theorem \ref{main-theorem}, as $n\tinf$, 
\[
 \sqrt{k} (\gamchabeta - \gamma_1) \limloi N(\lambda \mg,\sigg^2)
\]
where
\[
  \sigg^2=\displaystyle \frac{\gamma^2}{\pbet^2} \frac{p}{2 \pbet -1}  (1+\beta \gamma_1)^4 
 =\displaystyle \gamma_1^2\frac{p}{2p-1} (1+\beta \gamma_1)^2\frac{2p-1}{2 \pbet -1}
\]
and
\[  \mg = \left\{ \begar[c]{ll} 
  -\gamma^2 \beta_1 D_1 C^{-\gamma \beta_1} \pbet^{-1} (\pbet+\gamma \beta_1)^{-1} (1+\beta \gamma_1)^2 & \mbox{  if $\beta_1 \leq \beta_2 $},  \\
  0  & \mbox{   if $\beta_1 > \beta_2$.}
  \finar \right.
\]
\end{cor}

\begin{rmk}
Since  $\gamchapunk = \widehat{T}_k (0)= \widehat\gamma_{1,k}(0)$, taking $\beta =0$ in Theorem \ref{main-theorem} or in  Corollary \ref{corGammaW} entails the asymptotic normality for $ \gamchapunk$ when $p>1/2$, {\it i.e.} when $\gamma_2>\gamma_1$. When $\beta > 0$, the asymptotic normality for $\gamchabeta$ holds under the weaker assumption  $\pbet>1/2$, {\it i.e.} 
$\gamma_2>\gamma_1/(1+2\gamma_1\beta)$, and therefore allowing  for stronger censoring in the tail. On the other hand the restriction becomes worse for negative $\beta$.
\smallskip\\
 When $\beta_1 \leq \beta_2$ the absolute value of the asymptotic bias of $\gamchabeta$ is increasing in $\beta$.
 For a bias comparison for the case $\beta_1 > \beta_2$ one needs third order assumptions. On the other hand the asymptotic variance of $\gamchabeta$ is decreasing in $\beta$ as long as $p_\beta <1$ and is increasing as $p_\beta >1$. 
It is difficult to say anything in general about the comparison of the asymptotic mean-squared error of $\gamchabeta$ with respect to $\gamchapunk$. It is of course, when $\beta >0$ and  $p$ gets close to the value $1/2$,  in favor of $\gamchabeta$, at least from a theoretical point of view. 
\end{rmk}

\begin{rmk} From \citet{Einmahl08} it follows that the asymptotic variance of $\gamchapHillunk $ is given by ${1 \over k}{\gamma_1^2 \over p}$, which, for all $1/2<p<1$ is lower  than the asymptotic variance ${1 \over k}\frac{p \gamma^2_1}{2p-1}$  of $\gamchapunk$.
\smallskip\\
On the other hand, in case $\beta_1 \leq \beta_2$ it follows from \citet{BeirlantBardoutsosWetGijbels16}  that the  absolute value of the asymptotic bias of $\gamchapHillunk $ equals $(k/n)^{\gamma\beta_*}|m_{\gamma_1,0}| \frac{1+\gamma_1\beta_1}{1+\gamma \beta_1}$, which is larger than $(k/n)^{\gamma\beta_*}|m_{\gamma_1,0}|$ stated in the above theorem.  
\end{rmk}

\begin{rmk} The asymptotic distribution of $\gamchapunk$ in case $p \leq {1 \over 2}$, and in general of $\gamchabeta$ in case $p_{\beta} \leq {1 \over 2}$,  is not known. The authors conjecture  that asymptotic normality still holds, however with a slower rate than $k^{-1/2}$, presumably $k^{-p}$ when $p<1/2$, but the method of proof outlined below could not be carried through in that case.
\end{rmk}

Combining the asymptotic developments of $\gamchapunk$ and $\TB$ for $\beta=\beta_1$, which  are both weighted sums of the same i.i.d. random variables $p(E_i-1)-(\bI_{U_i\leq p}-p)$, and relying on the two-dimensional Lyapunov's CLT and the delta-method, it is now possible to deduce the following asymptotic normality result for the bias-reduced version of $\gamchapunk$ introduced in \citet{BeirlantMaribeVester18}. The proof is omitted for brevity.

\begin{cor} \label{corGammaBR}
Under the conditions of Theorem \ref{main-theorem} and assuming that $p>1/2$,  as $n\tinf$, we have 
\[
 \sqrt{k} (\gamBR - \gamma_1) \limloi N(0,\sigBR^2)
\]
where, with $\delta=p_{\beta_1}-p = \gamma\beta_1$, 
\[
  \sigBR^2 := \gamma_1^2 \frac{p}{2p-1} \frac{ (p+\delta)^2 ((p+\delta)^2+(1-p)^2+\delta+\delta^2)  } { \delta^2 (2p-1+\delta)(2p-1+2\delta) }.
\]
\end{cor}

\begin{rmk} While the asymptotic bias of $\gamBR$ is always $0$, its asymptotic variance is in  general  larger than those of the competing estimators.
\end{rmk}

\section{Proof of Theorem \ref{main-theorem}}  \label{Proofs}

 Let us introduce the following important notations with $1\leq i,j \leq k$:
\begeq{notationsksijuj}
 \ksi_j = \ksijdev 
 \makebox[1.2cm][c]{and} 
 \ui = \frac {i}{k+1}, 
\fineq
as well as the ratios
\begeq{notationsRF}
 \RFchapj = \RapFchapj 
 \makebox[1.2cm][c]{and} 
 \RFj = \RapFj .
\fineq
If we also define $ \ksijkbet  =  \ksi_j$ if $\beta =0$ and otherwise 
\[
 \ksijkbet  \; =  \; j  \left( k_{-\beta}\left(\frac{\Znmjun}{\Znmk}\right) -  k_{-\beta}\left(\frac{\Znmj}{\Znmk}\right)  \right)
  \; = \;
  \frac{j}{\beta} \left( \left(\frac{\Znmj}{\Znmk}\right)^{-\beta}  - \left(\frac{\Znmjun}{\Znmk}\right)^{-\beta} \right) 
\]
then, from $(\ref{defTB})$, we have 
\[
 \TB := \summ{j=2}{k} \RapFchapj  \frac{\ksijkbet}{j}
\]
where, using a Taylor expansion (of order $2$) , 
\begeq{decompXijkbeta}  
 \ksijkbet  \; =  \;
 \frac{j}{\beta} \left( \exp^{-\beta \log \left(\frac{\Znmj}{\Znmk} \right)} -  \exp^{-\beta \log \left(\frac{\Znmjun}{\Znmk}\right)} \right) 
  \; =  \;   
  \ksi_j    \left( \frac{\Znmjun}{\Znmk} \right)^{-\beta} + \beta \frac{\ksi^2_j}{2j}   \left( \frac{\tilde{Z}_{j,n}}{\Znmk} \right)^{-\beta},  
\fineq
for some variables $\tilde{Z}_{j,n}$ satisfying $\Znmj \leq \tilde{Z}_{j,n} \leq \Znmjun$.
\smallskip

The overall objective is to appropriately use the relation between the variables $\ksi_j$ and standard exponential order statistics  $\Ejn$ defined below, as well as between the ratios $\RFj$ and $(\Znmjun/\Znmk)^{-\beta}$ and uniform order statistics $\Vjk$ (with mean $\uj$) also defined below, in order to prove Theorem \ref{main-theorem}. Indeed, let $(Y_i)$ denote i.i.d. standard Pareto rv's defined by $Z_i=U_H(Y_i)$, and let  
\begeq{defEinVjk}
 \Ytildekmjun = Y_{n-j+1,n}/Y_{n-k,n}, \hspace{0.4cm} 
 \Vjk = 1/\Ytildekmjun,   \makebox[1.2cm][c]{ and }  
 \Ejn = j\log(Y_{n-j+1,n}/Y_{n-j,n}), \; 1\leq j\leq k.
\fineq
It is then known  that $(V_{1,k},\ldots,\Vjk, \ldots,V_{k,k})$ follows the distribution of the vector of order statistics of a standard uniform random sample of size $k$, and that the  variables $(E^{(n)}_1,\ldots,E^{(n)}_k)$ are jointly equal in distribution to a sample of size $k$ of  independent  standard exponential rv's.  
\smallskip

\citet{Beirlant02} showed that the rv's $ \ksi_j$ and $\Ejn$ are related as follows:
\begeq{decXij}
 \ksi_j  =  \ksi'_j  + \Rnj, \mbox{ where we define  }  \ksi'_j =(\gamma + \uj^{\gbetaet} \bnk)\Ejn , 
\fineq
where $\bnk$  is  asymptotically equivalent to $-\gamma^2 \beta_* D_* C^{-\gamma\beta_*} \left(  \frac{k+1}{n+1} \right)^{\gamma\beta_*}$, as $k,n \rightarrow \infty$ and $k/n \to 0$. Properties of the remainder  term $\Rnj$ will be detailed in Subsection \ref{sub-prelimlemmas}  . Equation $(\ref{decXij})$ thus implies that  
\begeq{decompXijkbetabis}
\ksijkbet = \ksijkbetp + \Rnjbet, 
\fineq
where 
\begin{eqnarray}  
\ksijkbetp & = &  \ksi'_j   \left( \frac{\Znmjun}{\Znmk} \right)^{-\beta} \label{defxipjkbeta}   \\
 \Rnjbet &= &  \Rnj  \left( \frac{\Znmjun}{\Znmk} \right)^{-\beta} +  \beta \frac{\ksi^2_j}{2j}   \left( \frac{\tilde{Z}_{j,n}}{\Znmk} \right)^{-\beta}.  \label{defRnjbeta} 
\end{eqnarray}  
\zun 
We can now start  breaking down $\TB - \frac{\gamma_1}{1+\gamma_1 \beta} $ into several terms  by writing:
\begin{align}
\TB - \frac{\gamma_1}{1+\gamma_1 \beta}
  & = \sumjdeuxk \RFchapj\frac{\ksijkbet}{j} - \frac{\gamma_1}{1+\gamma_1 \beta}   = \left( \sumjdeuxk \RFchapj  \frac{\ksijkbetp}{j} - \frac{\gamma_1}{1+\gamma_1 \beta} \right)  +    \sumjdeuxk \RFchapj \frac{\Rnjbet}{j}   \nonumber  \\ 
  & = \sumjdeuxk \left(\frac{\RFchapj}{\RFj} - 1\right) \RFj \frac{\ksijkbetp}{j}  \nonumber \\
  & \;\;\; + \ \left( \, \sumjdeuxk \RFj\frac{\ksijkbetp}{j} 
   \; - \; \frac{\gamma}{k+1}\sumjdeuxk \uj^{\pbet-1} \right)  
   \nonumber \\
   & \;\;\; + \ \left( \frac{\gamma}{k+1}\sumjdeuxk \uj^{\pbet-1} - \frac{\gamma}{\pbet} \right) \nonumber \\
   & \;\;\; +   \sumjdeuxk \RFchapj \frac{\Rnjbet}{j}  
  \nonumber \\
  & = \Tkn{1} + \Tkn{2} + R_n^{(0)}  + R_n^{(1)}  , \label{Tkn123etRn0}
\end{align}
with
\begin{align}
  \Tkn{1} 
  & =  \sumjdeuxk \left( \log\RFchapj - \log\RFj \right) \RFj \frac{\ksijkbetp}{j} \ + \ \sumjdeuxk \left\{ -\log \frac{\RFchapj}{\RFj} - \left( 1-\frac{\RFchapj}{\RFj} \right) \right\} \RFj\frac{\ksijkbetp}{j}  \nonumber  \\ 
  & =  \Tkn{1,1} \; + \; \Tkn{1,2}. 
 \label{decompTkn1}
\end{align}
The term $\Tkn{1,1}$ is introduced in order to make logarithms of  the Kaplan-Meier  product appear, leading to manageable sums.  Indeed, by definition of $\Fbarn$ we find that 
\[
  \log\RFchapj = \sum_{i=j}^k  \deltanmiun \log{\textstyle\left(\frac{i-1}{i}\right)} 
  \makebox[1.3cm][c]{ and } 
  \log\RFj= -\frac{1}{\gamma_1} \sum_{i=j}^k  \frac{\ksi_i}{i}  + \left(  \log\RFj+ \frac{1}{\gamma_1} \log \frac{\Znmjun}{\Znmk} \right) . 
  \]
 Consequently, defining the following important notations 
\begeq{defRFjbeta}
 \RFjbet  =   \RFj   \left( \frac{\Znmjun}{\Znmk} \right)^{-\beta}  \; i=2,\ldots,k,
\fineq
 and 
\begeq{defSikbeta}
  \Sikbet = \frac 1 i \summ{j=2}{i} \RFjbet \ksipjsurj  , \; i=2,\ldots,k,
\fineq
by inverting sums we obtain 
\begin{align}
 \Tkn{1,1} = \
 &  
 \sumideuxk \left[ \frac 1{\gamma_1} (\ksi_i - \gamma) + \left( \deltanmiun \,i\log{\textstyle\left(\frac{i-1}{i}\right)} \, + \, p \right) \right] \Sikbet -  \sumjdeuxk  \left(   \log\RFj+ \frac{1}{\gamma_1} \log \frac{\Znmjun}{\Znmk}   \right) \RFjbet \ksipjsurj  \nonumber \\
 = \ &   \Tkn{1,1,1}  \; - \; \Tkn{1,1,2}  .
 \label{decompTkn11}
\end{align} 

\noindent To summarize,  
\begin{equation}
\TB - \frac{\gamma_1}{1+\gamma_1 \beta}  = \Tkn{1,1,1}  - \Tkn{1,1,2}  + \Tkn{1,2} + \Tkn{2} + R_n^{(0)}  + R_n^{(1)}. \label{mastersum}
\zcinq
\end{equation}
Introducing now the additional notations 
\[
  c_i= 1+i\log \frac{i-1}{i},  \hspace{0.4cm}
  A_{i,n}=p(\Ein-1)- (\deltanmiun-p),  
  \makebox[1.3cm][c]{ and } 
  B_{i,n}=\frac{1}{\gamma_1} \bnk \ui^{\beta_* \gamma} \Ein,
\]
and using $(\ref{decXij})$,  one readily obtains the following formula for the main term $\Tkn{1,1,1}$ : 
\begeq{decTkn111}
\Tkn{1,1,1} = \sum_{i=2}^k A_{i,n} \Sikbet +  \sum_{i=2}^k B_{i,n} \Sikbet +  \sum_{i=2}^k \deltanmiun  c_i \Sikbet + \frac{1}{\gamma_1} \sum_{i=2}^k R_{n,i} \Sikbet . 
\fineq

\vspace{0.3cm}
 In the sequel, we will show  that the variables $\Sikbet$ can be  approximated appropriately by $\frac{\gamma}{\pbet}\frac{1}{k+1} \ui^{\pbet-1}$. 
 Also, as  it is explained in \citet{Einmahl08},  on one hand the parameter  $p =\gamma / \gamma_1=\frac{\gamma_2}{\gamma_1+\gamma_2}$ is the limit of $p(z)=\bP(\delta=1|Z=z)$ as $z\tinf$, and on the other hand  the original observations $(Z_i,\delta_i)_{i\leq n}$ have the same distribution  as the variables $(Z'_i,\delta'_i)_{i\leq n}$, where $(Z'_i)_{i\leq n}$ is an independent copy of the sequence  $(Z_i)_{i\leq n}$, \ $\delta'_i = \bI_{U_i \leq p(Z'_i)}$ and $(U_i)_{i\leq n}$ denotes some given i.i.d. sequence of standard uniform random variables (shortened to rv's), which are independent of the sequence $(Z'_i)_{i\leq n}$. We thus carry on the proof by considering from now on that the observations $\delta_i$ and $Z_i$ are related by the formula
  \[
   \delta_i=\bI_{U_i\leq p(Z_i)}.
  \]
Mimicking what is done in \citet{Einmahl08}, we  will later (see proof of Lemma  \ref{lemMartingalesEdelta}) approximate the rv's $\delta_{n-i+1,n}$ by  i.i.d Bernoulli rv's $\bI_{U_i \leq p}$.  
 
 \vspace{0.3cm}
 The main goal will thus be to  prove that the  term $ \sum_{i=2}^k A_{i,n} \Sikbet $ above is  (up to a bias term) close to the main  random term  appearing in Theorem \ref{main-theorem}
\begeq{mainterm}
\frac{\gamma}{\pbet}\frac{1}{k+1} \sumideuxk \left\{ p (E_i -1) - (\bI_{U_i \leq p}-p) \right\}  \ui^{\pbet-1}
\zdeux
\fineq
The other terms in $(\ref{decTkn111})$ will be bias or remainder terms, noting that the coefficients $c_i$ are close to $0$.
\zun\\
The second term $\Tkn{1,1,2}$ in $(\ref{mastersum})$ turns out to be adding to the bias since it only involves the slowly varying function $l_F$.
The treatment of the third term $\Tkn{1,2}$ above is very important since it strongly participates  to the approximation of a ratio of the form $\Fbarn(x)/\Fbarn(y)$ by the ratio $\Fbar(x)/\Fbar(y)$, for very large values of $x$ and $y$. Such approximation is  delicate. Invoking results from survival analysis, we will show however that $\Tkn{1,2}$ is a remainder term. 
\\\\
Next,   $\Tkn{2}$ is decomposed using the variables $(\Vjk)$ introduced in $(\ref{defEinVjk})$: 
\begeq{decompTkn2}
\Tkn{2} = \frac{\gamma}{k+1} \sumjdeuxk \left( \Vjkpbet -\ujpbet \right) \uj^{-1} + \sumjdeuxk   \Vjkpbet \frac{\ksi'_j -\gamma}{j} +  \sumjdeuxk \left( \RFjbet - \Vjkpbet \right) \ksipjsurj . 
 \fineq
 According to the definition of $\ksi'_j$, we can  see that the second term of this decomposition is close to 
 $  \frac{\gamma}{k+1} \sumjdeuxk (E_j -1)  \ujpmunbet $.
 While this is part of the main term described in $(\ref{mainterm})$, we will find in Proposition \ref{prop-Tkn2Tkn3-mild} that this term is neutralized by another part of $\Tkn{2}$, so  that $\Tkn{2}$ is just a bias term. Finally $R_n^{(0)}$ and $R_n^{(1)}$ will also turn out to be remainder terms. 
\medskip

The rest of the section is organised as follows. In subsection \ref{sub-prelimlemmas}, we set additional notations and state some preliminary approximation results  needed in the sequel. In subsection \ref{sub-lightcase} we state the asymptotic results for all terms in \eqref{mastersum} and conclude the proof.

\subsection{Additional  notations and important preliminary results}  \label{sub-prelimlemmas} 

\vspace{0.3cm}
\noindent
$\bullet$ First, in the sequel we will regularly work under the following event, for some $\alpha>1$ arbitrary close to $1$,
\begeq{defEnbeta}
 \calEn = \left\{ \; \forall 1\leq j\leq k \, , \ \alpha^{-1}\uj \leq \Vjk \leq \alpha \uj \; \right\},
\fineq
where $\uj$ and $\Vjk$ are defined in $(\ref{notationsksijuj})$ and $(\ref{defEinVjk})$.  According to \citet{ShorackWellner86} (chapter 8), for every $\alpha>1$ we have $\lim_{n\tinf}\bP(\calEn)=1$.  In the proof section, working "on the event $ \calEn$" will thus mean stating bounds  or results which are valid with an arbitrary large probability. 
\\\\
$\bullet$
Secondly,  the remainder term  $\Rnj$ defined in the second-order exponential representation of the log-spacings  $(\ref{decXij})$ satisfy, according to  Theorem 2.1  in  \citet{Beirlant02}, 
\begeq{proprieteRnj}
 \textstyle \left| \sum_{j=i}^k \frac{\Rnj}{j} \right| = o_{\bP}(\bnk \log_+ (\frac{1}{\ui})).
\fineq 
\\\\\
$\bullet$
Thirdly, under assumptions $(\ref{condFbar})$ and $(\ref{condGbar})$, since $Z_i=U_H(Y_i)$, one can show using $(\ref{condFbar})$ and $(\ref{condUH})$ that
\begeq{decRFj}
 \RFjbet = \RapFj  \left( \frac{\Znmjun}{\Znmk} \right)^{-\beta} =  \Vjkpbet (1+C_{j,k,\beta}) ,
  \fineq
where $C_{j,k,\beta}=Y_{n-k,n}^{-\gbetaet} D_{\beta}C^{-\gbetaet}(\Ytildekmjunexpo{-\gbetaet}-1)(1+\oPdeun)$ 
and $D_{\beta}=D - \gamma \beta D_*$ with $D=-\frac{\gamma}{\gamma_1}D_*$ if $\beta_2<\beta_1$, \ $D=D_1-\frac{\gamma}{\gamma_1}D_*$ if $\beta_1\leq \beta_2$. 
\\\\
$\bullet$
Finally, using R\'enyi representation (see for example $(4.3)$ in \citet{Beirlant04})  and  a Taylor expansion, one obtains  that for every $2\leq j\leq k$,  
\begeq{decompVjkpujkp}
 \Vjkpbet - \ujpbet \ = \ - \pbet \ujpbet \left( \summ{i=j}{k} \frac{\Ein-1}{i} \right) \; - \; \pbet \ujpbet\left( \summ{i=j}{k} \frac{1}{i} - \log\left(\frac{k+1}{j}\right) \right) 
 \; + \; \frac{\pbet^2}{2} \tilde V^{\pbet}_{j,k} \left(\log(\Vjk/\uj) \right)^2,
\fineq
where $\tilde V_{j,k}$ lies between $\Vjk$ and $\uj$. The combination of  $(\ref{decRFj})$ and $(\ref{decompVjkpujkp})$ 
thus means that the ratio $\RFjbet$ will be appropriately approximated by the deterministic weights $\ujpbet$.

\subsection{Asymptotics for the terms in \eqref{mastersum} and conclusion of the proof}  \label{sub-lightcase} 

The first result stated concerns the term $\Tkn{1,1,1}$, which contains the main term of the decomposition of $\TB - \frac{\gamma_1}{1+\gamma_1 \beta}$ (see relations $(\ref{decTkn111})$ and $(\ref{mainterm})$). 
\medskip
\begin{prop}
\label{prop-Tkn111}
Under the conditions of  Theorem \ref{main-theorem}, as $n \rightarrow  \infty$, we have
\begitem 
 \item[$(a)$]  $ \sqrt{k}\sum_{i=2}^k A_{i,n} \Sikbet   =  G_n + \lambda b_{\beta}  +o_{\bP} ( 1 )$,  \ \  
 where 
 $$ 
   b_{\beta}= - \gamma\frac{p}{\pbet} (1-p)  (D \gamma)_* \beta_* C^{-\gbetaet}/(\pbet+ \gamma \beta_*)   
   \makebox[1.3cm][c]{ and } 
   (D \gamma)_* =\left\{  \begar{lll}  \gamma_1 D_1  &  \mbox{ if }  & \beta_1 < \beta_2 \\
                                                - \gamma_2 D_2   & \mbox{ if }  & \beta_2 < \beta_1 \\
						\gamma_1 D_1 - \gamma_2 D_2 & \mbox{ if }  & \beta_1 = \beta_2 
						 \finar  \right.
 $$
 and   $G_n$ is equal in distribution to 
 $$
 \frac{\gamma}{\pbet}  \frac{1}{\sqrt{k}} \sum_{i=2}^k  \ui^{\pbet-1} \left( p(E_i-1)- (\bI_{U_i \leq p}-p) \right) ,
 $$
where $(E_i)$ and $(\delta_i)$ are independent iid samples with distributions standard exponential and standard uniform. The variable $G_n$ is asymptotically centred gaussian distributed with variance $ \sigT^2=\frac{\gamma^2}{\pbet^2} \frac{p}{2 \pbet -1}$. 
\item[$(b)$] $  \sqrt{k}\sum_{i=2}^k B_{i,n} \Sikbet  =     \lambda \bet  + o_{\bP} ( 1 )$ , \   where $\bet= - \gamma^2 \frac{p}{\pbet} D_* \beta_* C^{-\gbetaet}/(\pbet+ \gamma \beta_*)$. 
 \item[$(c)$]  $\  \sum_{i=2}^k \deltanmiun  c_i \Sikbet  =o_{\bP} ( k^{-1/2} )$ 
 \item[$(d)$]    $ \sum_{i=2}^k R_{n,i} \Sikbet  =o_{\bP} ( k^{-1/2} )$ 
\medskip
\finit
\end{prop}     

The following proposition concerns the terms  $R_n^{(0)}$, $R_n^{(1)}$, $\Tkn{1,2}$ and $\Tkn{1,1,2}$. The last two of these terms  result from the replacement of the ratios of Kaplan-Meier estimates $\RFchapj$ by the ratios of the true survival function values $\RFj$.
\begin{prop} 
\label{prop-Rn0-and-otherTterms}
Under the conditions of Theorem \ref{main-theorem},  as  $n \rightarrow  \infty$, 
\zun\\
\noindent $(a)$  $R_n^{(0)}  = o ( k^{-1/2}),$ \hspace{2.cm} 
$(b)$  $R_n^{(1)}  = o_{\bP} ( k^{-1/2}),$  \hspace{2.cm}
$(c)$  $\Tkn{1,2} = o_{\bP} ( k^{-1/2} ), $
\zun\\
 $(d)$   $\Tkn{1,1,2} =  D_1 (1+\oPdeun) \Znmkbetaun \sumjdeuxk \left(\left(\frac{\Znmjun}{\Znmk}\right)^{-\beta_1} - 1\right) \RFjbet \ksipjsurj 
 \; + \; \sumjdeuxk L_{n,j} \RFjbet \ksipjsurj$, where
$$
 0 \leq L_{n,j} \leq D_1^2 ( \Znmjunbetaun - \Znmkbetaun )^2 (1+\oPdeun).
$$
 Moreover,  $\Tkn{1,1,2} = b_{KM}   \left( k/n \right)^{\gbetaet} + o_{\bP} ( k^{-1/2} ) $,
where $ b_{KM}$ is equal to $-\frac{\gamma^2}{\pbet}  D_1 \beta_1 C^{-\gamma\beta_1} / (\pbet+ \gamma \beta_1)$  if $\beta_1 \leq \beta_2$ and  to $0$ if $\beta_1 > \beta_2$. 
\end{prop}

The last result concerns the behaviour of $\Tkn{2}$  : it turns out that  it only generates a bias term.
\begin{prop}
\label{prop-Tkn2Tkn3-mild}
We have 
\begin{align}
\label{lem-decompTkn2Tkn3}
   \Tkn{2} =
   & 
   - \frac{\pbet\,\gamma}{k+1} \sumjdeuxk (\Ejn-1)\left( \frac 1 j \sumideuxj \uipmunbet - \frac{1}{\pbet} \ujpmunbet \right) \; - \; 
   \frac{\pbet \,\gamma}{k+1} \sumjdeuxk  (\Ejn-1)  \ujpmunbet \left( \sum_{i=j}^k \frac{\Ein-1}{i} \right)  \nonumber  \\  
  & + \; \frac {\bnk} {k+1} \sumjdeuxk \uj^{\pbet-1+\gamma\beta_*}  \Ejn
  + \; \frac{\bnk}{k+1} \sumjdeuxk \frac 1{\uj} (\Vjkpbet-\ujpbet) \uj^{\gamma\beta_*}  \Ejn  \; + \; \sum_{j=2}^k \Vjkpbet C_{j,k,\beta} \ksipjsurj   \nonumber \\  
  &  - \; \frac{\pbet \,\gamma}{k+1} \sumjdeuxk \ujpmunbet \left( \sum_{i=j}^k \frac 1 i  - \log \frac{k+1}{j} \right)   \Ejn
   \; + \; \frac{(\pbet)^2\,\gamma}{2(k+1)} \sumjdeuxk \frac 1{\uj} \tilde{V}^{\pbet}_{j,k} \left(  \log\left(\frac{\Vjk}{\uj}\right) \right)^2 \Ejn .
\end{align}
Moreover, under the conditions of Theorem \ref{main-theorem}, when $n \rightarrow  \infty$ we have 
\[
      \Tkn{2}  =   \btildeet  \left( k/n \right)^{\gbetaet}  +o_{\bP} ( k^{-1/2} ), 
\]   
where $\btildeet  =  -  \frac{\gamma^2  \beta_* C^{-\gbetaet}}{\pbet+ \gamma \beta_*} (D_* + \frac{D_{\beta}}{\pbet}) $.
\end{prop}        
             
The proofs of all these results can be found in the Appendix. Now, since 
$$
 \sqrt{k}\left(\TB - \frac{\gamma_1}{1+\gamma_1 \beta} \right) = \sqrt{k}\Tkn{1,1,1} +\sqrt{k}\Tkn{1,2} - \sqrt{k}\Tkn{1,1,2} + \sqrt{k}\Tkn{1,1,3} + \sqrt{k}\Tkn{2} + \sqrt{k}R_n^{(0)} + \sqrt{k}R_n^{(1)}
$$
and assumption (\ref{conditionbiais}) holds, by combination of relation (\ref{decTkn111}) and propositions \ref{prop-Tkn111}, \ref{prop-Rn0-and-otherTterms}  and \ref{prop-Tkn2Tkn3-mild}, we have proved that Theorem \ref{main-theorem} holds, {\it i.e.} that 
$$
 \sqrt{k}\left(\TB - \frac{\gamma_1}{1+\gamma_1 \beta} \right) = G_n + \lambda \mT + o_{\bP} (1)  \stackrel{d}{\longrightarrow} N(\lambda \mT,\sigT^2),
$$
because it can be checked that $b_{\beta}+\bet-b_{KM}+\btildeet$ is actually equal to the value $\mT$ described in the statement of Theorem \ref{main-theorem}.
 

\zdeux \zdeux  
\section{Finite sample comparisons}  \label{FiniteSample}


In this section, we consider a comparison  (using finite sample simulations)  in terms of  observed bias and mean squared error (MSE) of the estimators considered in this paper : $\gamchapHillunk$, $\gamchapunk= \widehat\gamma_{1,k}(0)$, $\gamchabeta$ with $\beta \neq 0$,  and $\gamBR$.  For $\gamchabeta$, we consider three different values of $\beta$ ($-1$, $0.5$ and $1.5$). In the expression of $\gamBR$, the second order parameter $\beta_1$ of $F$ should be estimated. Instead, we proceed as in \citet{BeirlantMaribeVester18}  (see equations $(13)$ and $(14)$ therein) by  reparametrizing   $\beta_1 \gamchapunk$  by $-\rho_1$ and we consider  two different  values of  $\rho_1$ ($-1.5$ and $-2$) in the following formula 
$$
 \gamBR(\rho_1) = \gamchapunk   -  \frac{(1-\rho_1)^2 (1-2\rho_1)}{\rho_1^2}
\left(\widehat{T}_{k}\big(-\rho_1/\gamchapunk\big) -\frac{\gamchapunk}{1-\rho_1}\right).
$$
For the study of the sensitivity of this definition of $\gamBR(\rho_1)$ with respect to the choice of $\rho_1$, we refer to \citet{BeirlantMaribeVester18}.
\zdeux

We consider  two classes of heavy-tailed distributions for the target and censoring variables $X$ and $C$ : 
\begitem 
\item Burr$(\theta,\beta,\lambda)$  with d.f. $1-(\frac{\theta}{\theta+x^{\beta}})^{\lambda}$, which extreme value index  is $\frac{1}{\lambda \beta}$. 
\item Fr\'echet$(\gamma)$  with d.f. $\exp(-x^{-1/\gamma})$, which extreme value index is $\gamma$. 
\finit

For each considered distribution, $2000$ random samples of length $n=500$ were generated ; median bias and MSE of the  above-mentioned estimators are plotted against different values of $k_n$, the number of excesses used. 
\zun

We considered two cases : a Burr distribution censored by another Burr distribution (Fig.1), a Fr\'echet distribution censored by another Fr\'echet distribution (Fig.2).
In each case, we considered a situation with $p>1/2$, which corresponds to weak censoring in the tail,  and the reverse situation with $p<1/2$, which corresponds to strong censoring. 
In  the Burr case,  we also considered situations with $\beta_1 < \beta_2$,  and  reverse situations with $\beta_1> \beta_2$.  Indeed, for Fr\'echet distibutions,  $\beta_1$ is always larger that $\beta_2$ in the case $p>1/2$ and $\beta_1$ is always lower that $\beta_2$ in the case $p<1/2$.
\zun

\begin{figure}[htp]
    \centering
    \subfigure[Burr$(10,2,5)$ censored  by Burr$(10,4,1)$]{
          \label{BurrBurrfaible}
          \includegraphics[height=3.5cm,width=.46\textwidth]{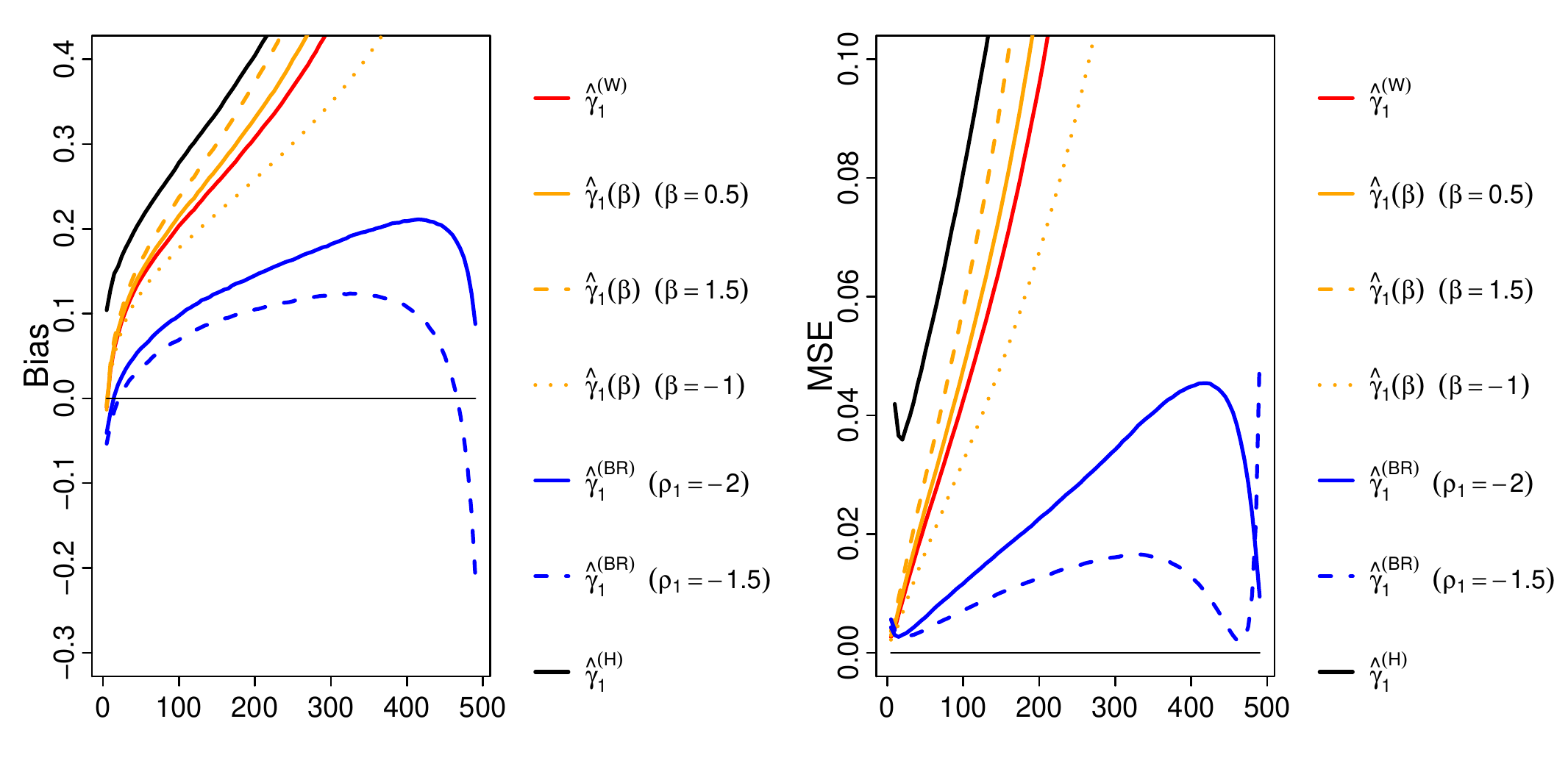}}
     \hspace{.1in}
     \subfigure[Burr$(10,2,2)$ censored  by Burr$(10,5,2)$]{
          \label{BurrBurrforte}
          \includegraphics[height=3.5cm,width=.46\textwidth]{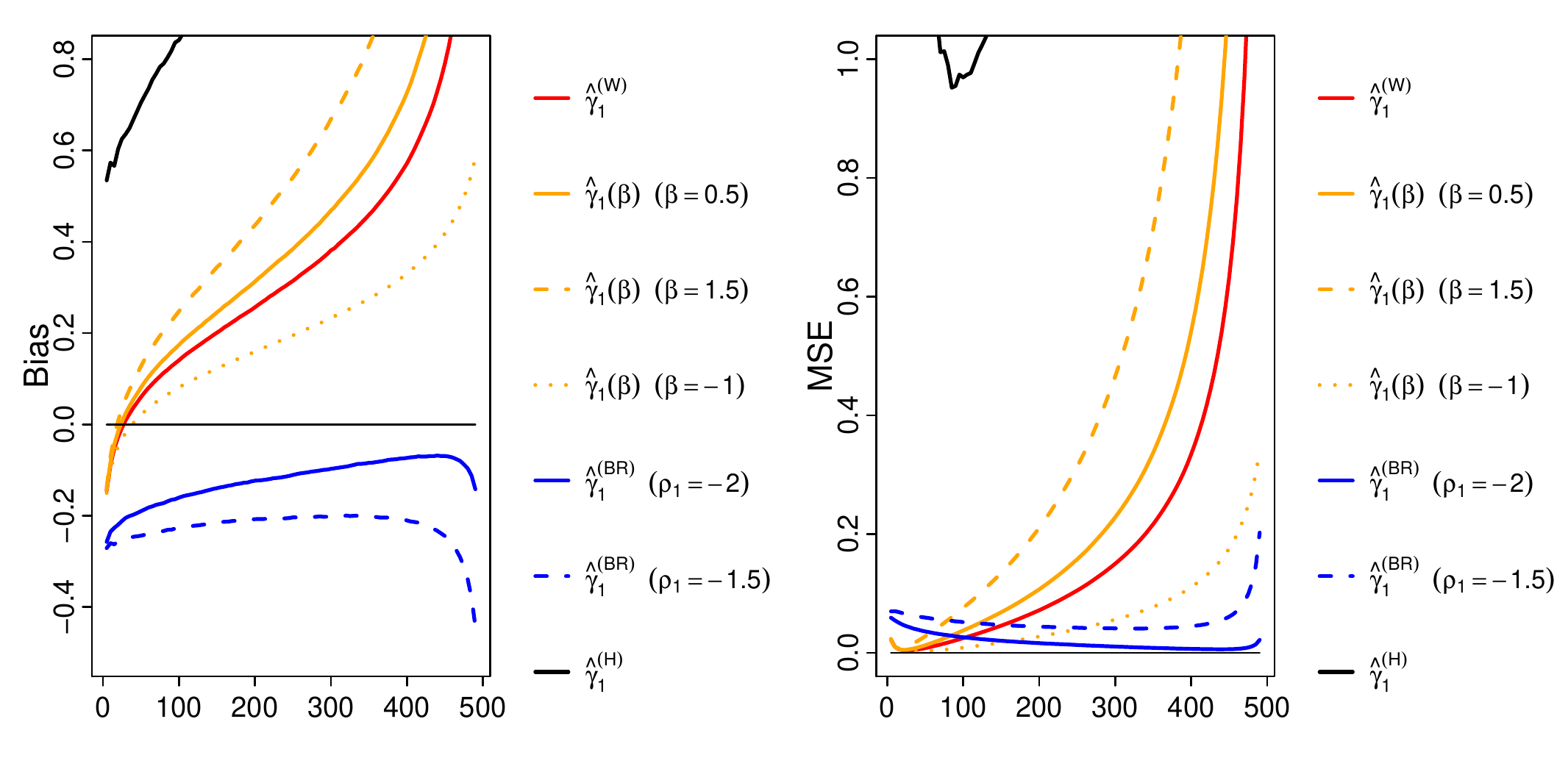}}
       \subfigure[Burr$(10,5,2)$ censored  by Burr$(10,2,2)$]{
          \label{BurrBurrfaible2}
          \includegraphics[height=3.5cm,width=.46\textwidth]{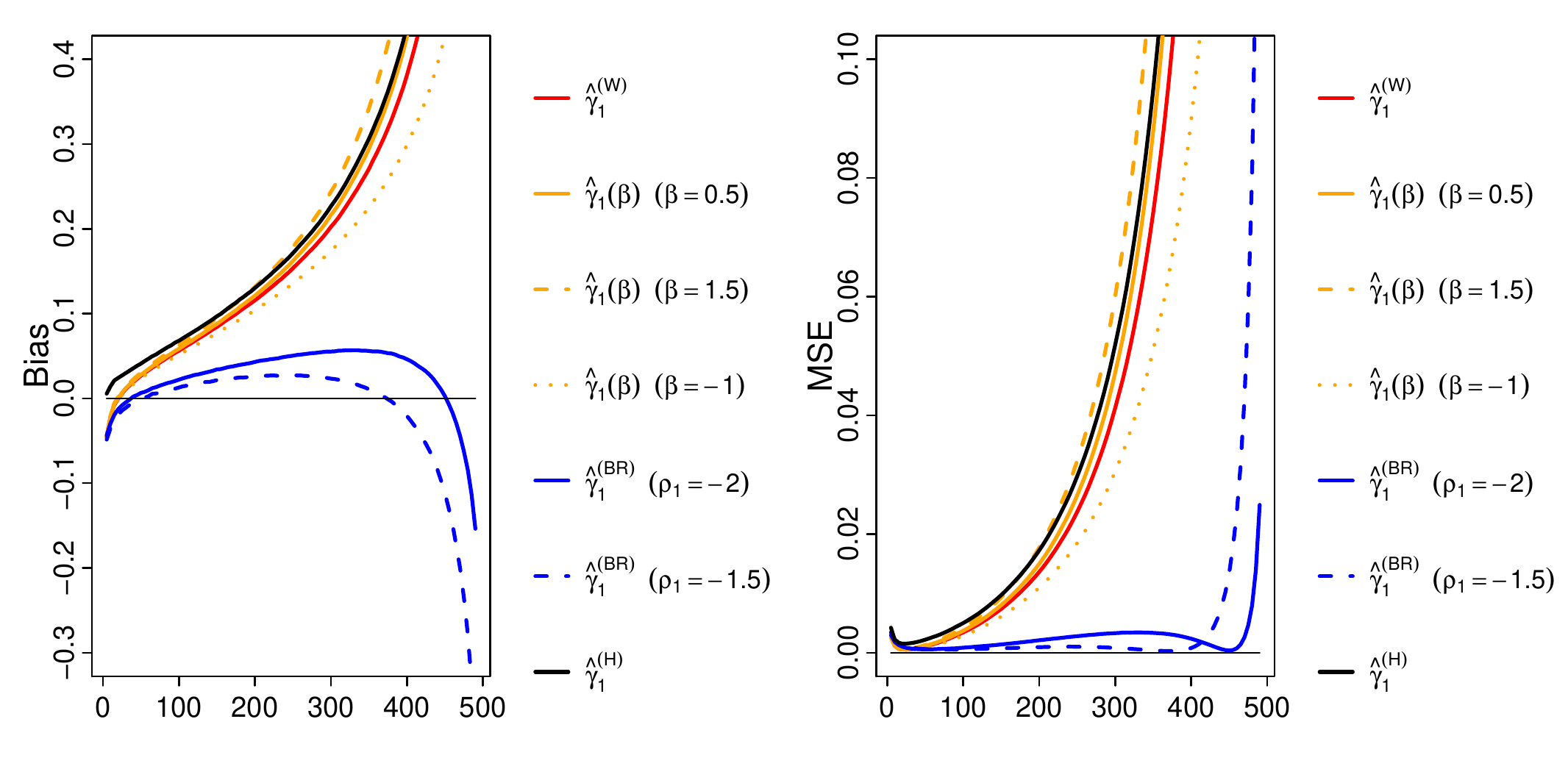}}
     \hspace{.1in}
     \subfigure[Burr$(10,4,1)$ censored  by Burr$(10,2,5)$]{
          \label{BurrBurrforte2}
          \includegraphics[height=3.5cm,width=.46\textwidth]{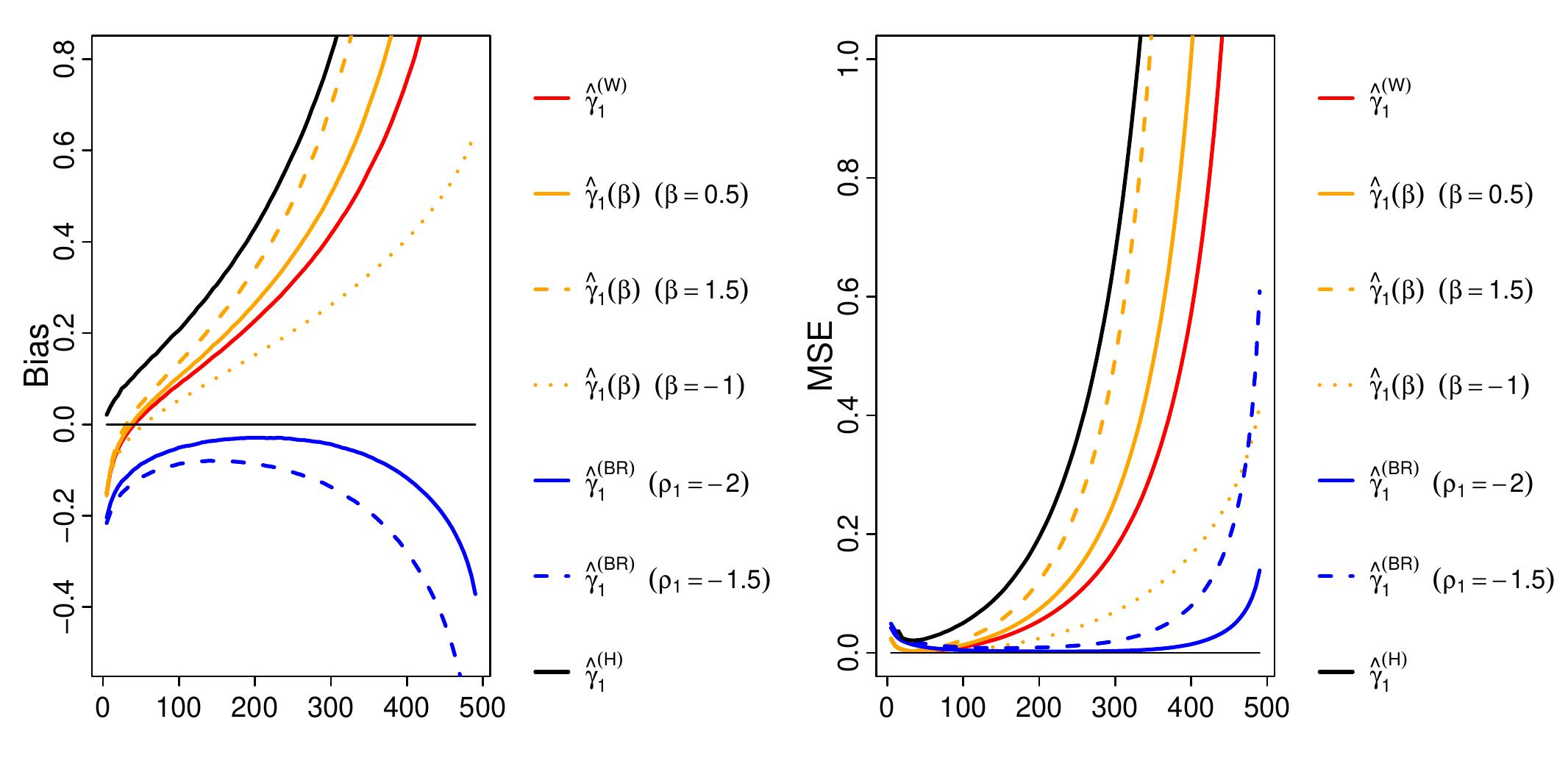}}
    \caption{Comparison of bias and MSE for $\gamchapHillunk$, $\gamchapunk= \widehat\gamma_{1,k}(0)$, $\gamchabeta$  and $\gamBR(\rho_1)$ for a Burr distribution censored by another Burr distribution : (a) $\beta_1 =2 < \beta_2=4$ and $p>1/2$, \ (b) $\beta_1 =2 < \beta_2=5$ and $p<1/2$, \ (c) $\beta_1 =5 > \beta_2=2$ and $p>1/2$, \ (d) $\beta_1 =4 > \beta_2=2$ and $p<1/2$ }.
    \label{burrburr}
\end{figure}

\begin{figure}[htp]
    \centering
    \subfigure[Fr\'echet$(1/4)$ censored  by Fr\'echet$(1/2)$]{
          \label{FFfaible}
          \includegraphics[height=3.5cm,width=.46\textwidth]{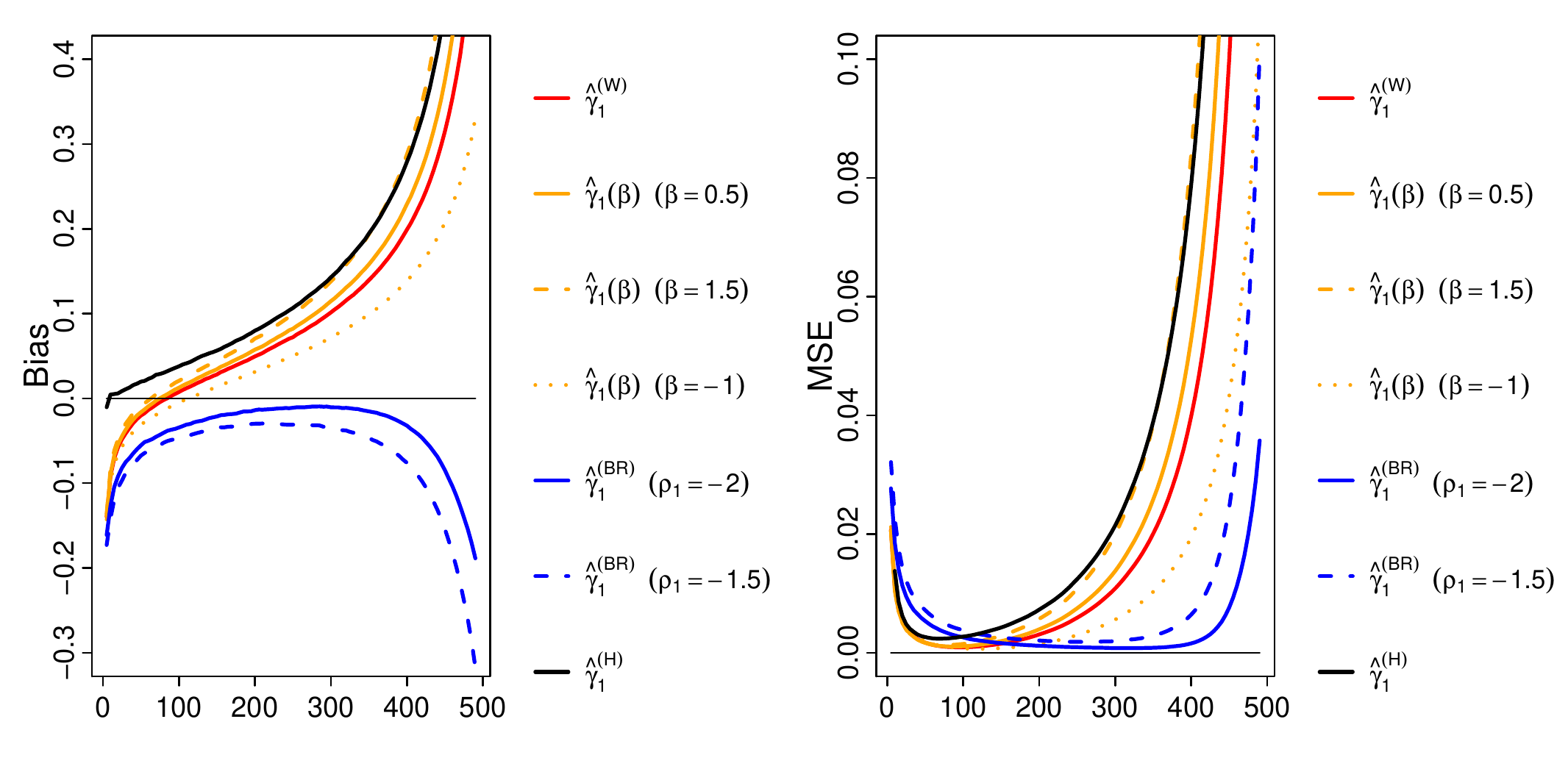}}
     \hspace{.1in}
     \subfigure[Fr\'echet$(1/2)$ censored  by Fr\'echet$(1/4)$]{
          \label{FFforte}
          \includegraphics[height=3.5cm,width=.46\textwidth]{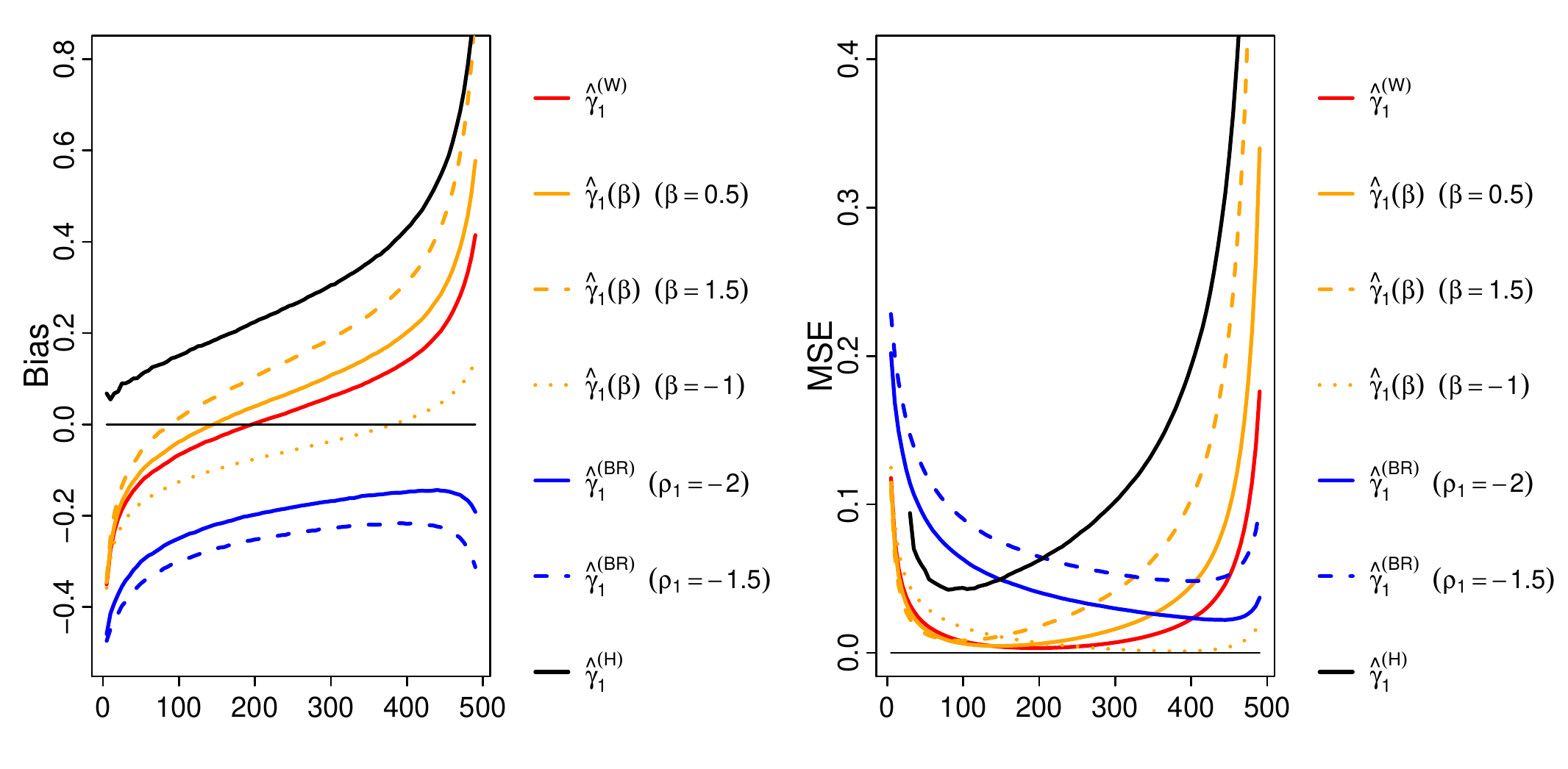}}
   \caption{Comparison of bias and MSE for $\gamchapHillunk$, $\gamchapunk= \widehat\gamma_{1,k}(0)$, $\gamchabeta$  and $\gamBR(\rho_1)$   for a Fr\'echet  distribution censored by another Fr\'echet  distribution : (a) $\beta_1 =4 > \beta_2=2$ and $p>1/2$, \ (b) $\beta_1 =2 < \beta_2=4$ and $p <1/2$. }
     \label{FF}
\end{figure}

This small simulation study shows  that the MSE of $\gamchabeta$ is globally decreasing with lower values of $\beta$, even when the condition $p_{\beta} > {1 \over 2}$ for the above asymptotic normality result is not met, as in the case with $\beta =-1$ and $p < {1 \over 2}$. This is probably due to the decreasing bias with decreasing $\beta$, the bias being the dominating component in the MSE.

On the other hand  $\gamBR$ overall reduces the MSE  for most $k$,  except in the heavy censoring Fr\'echet case.  The non-optimal behavior for small  values of $k$  is a well-known characteristic of bias reduced estimators. In \citet{BeirlantMaribeVester18} a penalized bias reduction technique was proposed to remedy this fact.  

\zun\zdeux\zdeux

\newpage

\section{Appendix} 


\subsection{Useful Lemmas}

Some of the following ten Lemmas are used several times in the proof of Propositions \ref{prop-Tkn111},   \ref{prop-Rn0-and-otherTterms}   and  \ref{prop-Tkn2Tkn3-mild}. 

\begin{lem}
\label{lem-majorationsdeterministes}
For any integer $i  \geq 2$ and every $k\geq i$ , we have
\begeq{ineqci} 
  c_i= 1 + i \log \frac{i-1}{i} \in \left[-\frac 1 i , 0\right]
\fineq
\begeq{ineq1surj}
  \sum_{j=i}^k \frac 1 j  - \log \frac{k+1}{i}  \in \left[0, \frac 1 i \right]
\fineq
Moreover, for any given $a\in ]0,1[$, there exist  some  positive constants $C_1<C_2$  such that, for all $2\leq i\leq k$
\begeq{ineqdik}
   d_{i,k}= \left( \frac 1 i \sum_{j=2}^i \uj^{-a}  - \frac{1}{1-a} \ui^{-a}  \right) \in \left[-\frac{C_2}{\ui (k+1)^{1-a}}, -\frac{C_1}{\ui (k+1)^{1-a}}\right],
\fineq
as well as, if $a<0$, 
\begeq{ineqdik2}
   d_{i,k}= \left( \frac 1 i \sum_{j=2}^i \uj^{-a}  - \frac{1}{1-a} \ui^{-a}  \right) \in \left[-\frac{1}{\ui (k+1)}, \frac{-a}{\ui (k+1)}\right].
\fineq
\end{lem}

\begin{lem} 
\label{lem-weightedsums}
For any $a<1$, we have,  as $n\tinf$,
\begeq{sumdet1}
 \frac 1 k \sum_{j=1}^k \uj^{-a}   \rightarrow \frac{1}{1-a}, 
\fineq
and, under assumptions $(\ref{condFbar})$ and $(\ref{condGbar})$, 
\begeq{sumxi}
 \frac 1 k \sum_{j=1}^k \uj^{-a} \  \ksi_j \limprob \frac{\gamma}{1-a},  
\fineq
(equation $(\ref{sumxi})$ also holds for $\ksi'_j$ instead of $\ksi_j$) and, if $X_j$ denotes either $E_j$, $E_j-1$ or $|E_j-1|$, where $(E_j)$ are  standard exponential iid random variables, then we have
\begeq{sumEi}
  \frac 1 k \sum_{j=1}^k \uj^{-a} \  X_j \limprob \frac{\bE(X_1)}{1-a}, \mbox{ as  } n \rightarrow +\infty.
\fineq
\end{lem}

\begin{lem} 
\label{lem-weightedsumsbis}
For any $a>1$, we have,  as $n\tinf$,
\begeq{sumdet}
 \sum_{j=1}^k j^{-a} \rightarrow  \zeta(a) \mbox{ \ \ as  } k \rightarrow +\infty,
 \fineq
 where $ \zeta$ is the Riemann Zeta function.  Moreover,  for any $\delta >0$, under  $(\ref{condFbar})$ and $(\ref{condGbar})$, 
\begeq{sumxibis}
\frac{1}{k^{a+\delta}} \sum_{j=1}^k \uj^{-a} \  \ksi_j \limprob 0,  \mbox{ as  } n \rightarrow +\infty,
\fineq
(equation $(\ref{sumxibis})$ also holds for $\ksi'_j$ instead of $\ksi_j$) and, if $(X_j)$ is  a sequence of i.i.d. random variables such that $\bE(|X_1|) <+\infty$,  then 
\begeq{sumEibis}
  \frac{1}{k^{a+\delta}} \sum_{j=1}^k \uj^{-a} \  X_j \limprob 0, \mbox{ as  } n \rightarrow +\infty.
\fineq
\end{lem}

\begin{lem}   \label{lemmaVjk}
If $(\Vjk)_{1\leq j \leq k}$ are the order statistics of $k$ standard uniform random variables then, for any $0<\delta<1$ and  $a>0$,  we have, as $k\tinf$,
\begeq{lemmaVjkpartie2}
 \sqrt{k} \max_{2\leq j\leq k} \frac{|\Vjk^a-\uj^a|}{\uj^{a-1/2-\delta/2}} = O_{\bP}(1).
\fineq
\end{lem}

\begin{lem}  \label{maxexpo}
If $(E_j)$ are  standard exponential iid random variables, then  
$ \max_{2\leq j\leq k}  |E_j| =O_{\bP}(\log k)$.
\end{lem}

\begin{lem} \label{PotterBounds} (See  \cite{HaanFerreira06}  Proposition B.1.9) \\
Suppose $f \in RV_{\alpha}$. If $x>0$ and $\delta_1,\delta_2  >0$ are given, then there exists $t_0=t_0(\delta_1,\delta_2)$ such that for any $t\geq t_0$ satisfying $tx \geq t_0$, we have
\[
(1-\delta_1) x^{\alpha} \min(x^{\delta_2}, x^{-\delta_2})  < \frac{f(tx)}{f(t)} < (1+\delta_1) x^{\alpha} \max(x^{\delta_2}, x^{-\delta_2}) .
\]
 If $x \geq  1$, then there exists $t_0=t_0(\epsilon)$ such that for every $t\geq t_0$, 
\begeq{bornesPotter}
(1-\epsilon) x^{\alpha-\epsilon} < \frac{f(tx)}{f(t)} < (1+\epsilon) x^{\alpha+\epsilon} .
\fineq
\end{lem}

\begin{lem} \label{lemMartingalesEE}
If $(E_i)_{i\leq k}$ are  standard exponential iid random variables, then if $\pbet> 1/2$, as $n\tinf$,
\begeq{EmunEmun}  
\frac{1}{\sqrt{k}} \sum_{i=3}^k (E_i-1) \ \left\{ \frac 1 i \sum_{j=2}^{i-1} \ujpmunbet  (E_j -1) \right\}  \limprob  0, 
\fineq
\begeq{EmunE}   
\frac{1}{k} \sum_{i=3}^k  (E_i -1) \ \left\{ \frac 1 i   \sum_{j=2}^{i-1} u_{j,n}^{\pbet+d-1} E_j \right\}  \limprob  0   \hs{0.5cm} (\mbox{for any $d\geq 0$}) 
\fineq
\begeq{EmunEbis}   
\frac{1}{k} \sum_{i=3}^k  \ui^{{\beta_* \gamma} } E_i  \ \left\{ \frac 1 i   \sum_{j=2}^{i-1} \ujpmunbet (E_j -1) \right\}  \limprob  0 . 
\fineq
\begeq{EmunEmunE}   
\frac{1}{\sqrt{k}} \sum_{i=4}^k  (E_i -1) \ \left\{ \frac 1 i   \sum_{l=3}^{i-1} (E_l-1) \left( \frac{1}{l} \sum_{j=2}^{l-1} u_{j,k}^{\pbet-1} E_j \right) \right\}  \limprob  0   
\fineq

\end{lem}

\begin{lem} \label{lemMartingalesEdelta}
With $\deltanmiun$ and $\Ejn$ being respectively defined in the introduction and in equation   $(\ref{defEinVjk})$,   if $\pbet> 1/2$ then, we have, under assumptions $(\ref{condFbar})$ and $(\ref{condGbar})$, as $n\tinf$,
\begeq{Edelta}  
\frac{1}{\sqrt{k}} \sum_{i=2}^k (\deltanmiun -p) \ \left\{ \frac 1 i \sum_{j=2}^{i-1} \ujpmunbet  (\Ejn-1)  \right\} \limprob  0 
\fineq
\begeq{deltaE}   
\frac{1}{k} \sum_{i=2}^k  (\deltanmiun -p) \ \left\{ \frac 1 i   \sum_{j=2}^{i-1} u_{j,n}^{\pbet+d-1} \Ejn  \right\} \limprob  0  \hs{0.5cm} (\mbox{for any $d\geq 0$}). 
\fineq
\begin{eqnarray}   
\frac{1}{\sqrt{k}} \sum_{i=3}^{k-1}   (\deltanmiun -p)  \left( \frac 1 i \sum_{j=2}^{i-1}  \ujpmunbet \Ejn \right) \left( \summ{l=i+1}{k}  \frac{\Eln-1}{l} \right) & = & \nonumber \\
\frac{1}{\sqrt{k}} \sum_{i=4}^k  (\Ein -1) \ \left\{ \frac 1 i  \sum_{l=3}^{i-1} (\delta_{n-l+1,n}-p) 
\left( \frac 1 l \sum_{j=2}^{l-1} \ujpmunbet \Ejn \right) \right\}  & \limprob &  0
\label{EmunEdeltamp}   
\end{eqnarray}
\begeq{deltaEmunE}   
\frac{1}{\sqrt{k}} \sum_{i=4}^k  (\deltanmiun -p) \ \left\{ \frac 1 i   \sum_{l=3}^{i-1} (\Eln-1) \left( \frac{1}{l} \sum_{j=2}^{l-1} u_{j,k}^{\pbet-1} \Ejn \right) \right\}  \limprob  0   
\fineq
\end{lem}

\begin{lem} \label{fonctionp} 
Let $p(z)=\bP(\delta=1|Z=z)$. Under the Hall model (conditions $(\ref{condFbar})$ and $(\ref{condGbar})$), 
\begeq{pHall} 
p\circ U_H(x)= p+p(1-p) (D \gamma)_* \beta_* C^{-\gamma \beta_*} x^{-\gamma \beta_*} (1+o(1)). 
\fineq 
Moreover, $(\ref{pHall})$ and $(\ref{conditionbiais})$ imply that 
\begeq{limp1} 
 \frac{1}{\sqrt{k}}  \sumideuxk \uipmunbet \left(  p\circ U_H(n/i) -p  \right)  \rightarrow   \lambda \alpha_{\beta}, 
\fineq 
where $\alpha_{\beta} = \frac{1}{\pbet+\gamma \beta_*} p(1-p) (D \gamma)_* \beta_* C^{-\gamma \beta_*}$. 
\end{lem}

\begin{lem} \label{Lemmetriplesomme} 
Using the notations introduced earlier, we have,  under assumptions $(\ref{condFbar})$ and $(\ref{condGbar})$ and if $\pbet> 1/2$, as $n\tinf$,
\[
\sqrt{k} \ \sumideuxk A_{i,n} \frac 1 i  \left( \sum_{j=2}^i  (\Vjk^{\pbet} -\uj^{\pbet}) \frac{\Ejn}{j} \right)   \limprob  0 .
\]
\end{lem}

\noindent 
We now prove one after the other the Propositions \ref{prop-Tkn111}, \ref{prop-Rn0-and-otherTterms} and \ref{prop-Tkn2Tkn3-mild}, then we will deal with the proofs of the different Lemmas in subsections \ref{subsectionproofLemme1} to \ref{subsectionProofLemmesMineurs}.

\subsection{Proof of Proposition \ref{prop-Tkn111}} \label{subsectionProp1}

\subsubsection{Proof of part $(a)$}
 \label{preuve-morceauprincipal-sommeAinSik}
 
This subsection is  devoted to the study of  $\sumideuxk A_{i,n}\Sikbet$, which we divide in three parts, using statement $(\ref{decRFj})$  : 
\[
 I_{1,n} + I_{2,n} + I_{3,n} = 
 \sumideuxk A_{i,n} \left( \frac 1 i  \sum_{j=2}^i  \uj^{\pbet} \ksipjsurj  \right)   +
 \sumideuxk A_{i,n}  \left( \frac 1 i \sum_{j=2}^i  (\Vjk^{\pbet} -\uj^{\pbet}) \ksipjsurj  \right) +
 \sumideuxk A_{i,n} \left(  \frac 1 i \sum_{j=2}^i  \Vjk^{\pbet} C_{j,k,\beta}\ksipjsurj  \right). 
\]
From $I_{1,n}$ will come the asymptotically gaussian part of $\sumideuxk A_{i,n}\Sikbet$, plus a bias term, and the other two $ I_{2,n}$ and $ I_{3,n}$ will be remainder terms. We will first give details about $I_{1,n}$, and then come back to $I_{2,n}$ and $I_{3,n}$ later. 
\zdeux
In order to deal with $I_{1,n} $, we begin by using relation $(\ref{decXij})$  to write $\ksi'_j $ as $\gamma + \gamma  (\Ejn -1)+ \uj^{\gbetaet} \bnk \Ejn$, which divides $I_{1,n}$ in  three different terms $I_{1,n}=I^{(1)}_{1,n}+I^{(2)}_{1,n}+I^{(3)}_{1,n}$. 
\zdeux

Our first task will be to deal with the main term of the theorem, $I^{(1)}_{1,n}$.  Recalling that $A_{i,n}=p(\Ein-1)- (\deltanmiun-p)$, where $ \delta_i=\bI_{U_i\leq p(Z_i)}$ with $(U_i)$ uniformly distributed and independent of $(Z_i)$ and $U_{n-i+1,n}$ denotes the uniform variable associated to $\deltanmiun$, this first term is equal to  
\begin{eqnarray*}
 I^{(1)}_{1,n} & = & \frac{\gamma}{k+1}   \sum_{i=2}^k A_{i,n}  \left(\frac{1}{i}  \sum_{j=2}^i \ujpmunbet  \right) 
 \\ 
 & =  &  
\frac{\gamma}{\pbet}  \frac{1}{k+1}  \sum_{i=2}^k  \uipmunbet \left( p(\Ein-1)- (\bI_{U_{n-i+1,n} \leq p}-p) \right)  
 \\
 & & \hspace*{1.cm}
 - \
\frac{\gamma}{\pbet}  \frac{1}{k+1}  \sum_{i=2}^k  \uipmunbet \left( \bI_{U_{n-i+1,n} \leq p(\Znmiun)} - \bI_{U_{n-i+1,n} \leq p} \right) 
 \ + \
 \frac{\gamma}{k+1}   \sum_{i=2}^k A_{i,n}  d_{i,k}
 \\
  & = & W_{k,n} + B_{k,n} + R_{k,n} 
\end{eqnarray*}
where we define $d_{i,k} = \frac 1 i  \sum_{j=2}^i \ujpmunbet -  \frac{1}{\pbet}\uipmunbet$.  To sum up what we have found so far,
\[
\sumideuxk A_{i,n}\Sikbet = (W_{k,n} + B_{k,n} + R_{k,n})\; + \; (I^{(2)}_{1,n}+I^{(3)}_{1,n}) \; + \; I_{2,n} \; + \; I_{3,n}.
\]

Introducing a sequence $(E_i)$ of independent standard exponential variables, independent of the sequence $(Z_i)$, we can write that 
\[
 W_{k,n} \equloi \frac{\gamma}{\pbet}  \frac{1}{k+1}  \sum_{i=2}^k  \uipmunbet \left( p(E_i-1)- (\bI_{U_i \leq p}-p) \right) \mbox{ \ and \ } 
 B_{k,n} \equloi -\frac{\gamma}{\pbet}  \frac{1}{k+1}   \sum_{i=2}^k  \uipmunbet \left( \bI_{U_i \leq p(\Znmiun)} - \bI_{U_i \leq p} \right), 
\]
We prove easily  that $\bV ar(\sqrt{k} W_{k,n})$ is equivalent to  the variance $\sigT^2$ defined in the statement of Theorem \ref{main-theorem}, and that, using Lyapunov's CLT,  we have $\sqrt{k} W_{k,n}  \limloi  N(0,\sigT^2)$.  
\zun 

Let us now deal with the term  $B_{k,n} \equloi B^{(1)}_{k,n} + B^{(2)}_{k,n} $, where 
\begin{eqnarray}  
B^{(1)}_{k,n}  & =  & - \frac{\gamma}{\pbet}  \frac{1}{k+1}  \sum_{i=2}^k  \uipmunbet \left( \bI_{U_i \leq p\circ U_H(Y_{n-i+1,n})} -  \bI_{U_i \leq p\circ U_H(n/i)}  \right) \nonumber \\
B^{(2)}_{k,n}  &= &   - \frac{\gamma}{\pbet}  \frac{1}{k+1}   \sum_{i=2}^k  \uipmunbet \left(  \bI_{U_i \leq p\circ U_H(n/i)} -  \bI_{U_i \leq p} \right). \nonumber 
\end{eqnarray} 
Following the method used for the treatment of the terms $T_{1,k}$ and $T_{2,k}$ in \citet{Einmahl08}, and using the LLN result found for instance in \citet{ChowTeicher1997} page 356, we can prove that $\sqrt{k} B^{(1)}_{k,n}   \limprob  0$ and that (using $(\ref{limp1})$, wherein constant $\alpha_{\beta}$ is defined) $\sqrt{k} B^{(2)}_{k,n}  \limprob - \frac{\gamma}{\pbet}  \lambda \alpha_{\beta}= \lambda  b_{\beta}$.
\zun

Concerning now the last term $R_{k,n}$ of $I^{(1)}_{1,n}$, if $\pbet<1$, according to inequality $(\ref{ineqdik})$ in Lemma \ref{lem-majorationsdeterministes}, there exists some constant $c>0$ such that 
\[
 \sqrt{k} |R_{k,n}| \leq 
 \sqrt{k} \frac{c \gamma}{(k+1)^{\pbet+1}}   \sum_{i=2}^k |A_{i,n}| \frac{1}{\ui} 
 \leq 
 O(1) k^{-(\pbet-1/2-\delta)} \frac 1 k  \sum_{i=2}^k |A_{i,n}|  \ui^{\delta-1}, 
\]
for a given $\delta >0$. But $ |A_{i,n}| \leq p |\Ein -1| + 1 \leq \Ein +2$, and therefore, taking $\delta$ small enough,  $\sqrt{k} R_{k,n} = o_{\bP}(1)$ according to properties $(\ref{sumEi})$ and $(\ref{sumdet1})$ (in Lemma \ref{lem-weightedsums}, with $a=1-\delta$)  and to the assumption $\pbet>1/2$. When $\pbet >1$, the treatment is similar, using $(\ref{ineqdik2})$ instead of $(\ref{ineqdik})$.  We have thus finished to prove that $\sqrt{k}I^{(1)}_{1,n}$ converges in distribution to $N(\lambda b_{\beta},\sigT^2)$. All the remaining terms in this subsection will now be proved to be negligible. 
\zdeux

Let us now consider the second term $I^{(2)}_{1,n}$ of $I_{1,n}$. Separating $j<i$ and $j=i$, we have  
\[
 I^{(2)}_{1,n} =
 \frac{\gamma}{k+1} \sum_{i=3}^k A_{i,n} \frac 1 i  \left( \sum_{j=2}^{i-1}  \uj^{\pbet-1} (\Ejn -1) \right) +   \frac{\gamma}{(k+1)^2} \sumideuxk A_{i,n}    \ui^{\pbet-2} (\Ein -1). 
\]
The first term  is shown to be $o_{\bP}(k_n^{-1/2})$ by separating $A_{i,n}$ in its $(\Ein-1)$ and $(\deltanmiun-p)$ parts and relying on properties $(\ref{EmunEmun})$ and $(\ref{Edelta})$ stated in Lemmas \ref{lemMartingalesEE} and \ref{lemMartingalesEdelta}. The second one is easy to handle using $(\ref{sumEibis})$ and $\pbet>1/2$ ; it is then omitted. 
\zdeux

Similarly, the third  term $I^{(3)}_{1,n}$ of $I_{1,n}$ is, again seperating  $j<i$ and $j=i$, 
\[
 I^{(3)}_{1,n} =
 \frac{\bnk}{k+1} \sum_{i=3}^k  A_{i,n} \frac 1 i  \left( \sum_{j=2}^{i-1}  \uj^{\pbet+\gamma \beta_*-1} \Ejn  \right) +   \frac{\bnk}{(k+1)^2} \sumideuxk A_{i,n}    \ui^{\pbet-2+\gamma \beta_*} \Ein . 
\]
Since $\sqrt{k} \bnk$ converges to a constant, the first term is $o_{\bP}(k_n^{-1/2})$ by  using properties $(\ref{EmunE})$ and $(\ref{deltaE})$ (with $d=\gamma\beta_*$) stated in Lemmas  \ref{lemMartingalesEE} and \ref{lemMartingalesEdelta}. Again,  the second one is easy to handle using $(\ref{sumEibis})$. 
\zdeux

Now that we have finished with $I_{1,n}$, we turn to the term $I_{2,n}$. The decomposition of $ \ksi'_j $ in  $(\ref{decXij})$ and the fact that $\sqrt{k} \bnk$ converges imply that 
\[
\sqrt{k} I_{2,n}  = \gamma \sqrt{k} \sumideuxk A_{i,n} \frac 1 i  \left( \sum_{j=2}^i  (\Vjk^{\pbet} -\uj^{\pbet}) \frac{\Ejn}{j} \right) + O(1) \sumideuxk A_{i,n} \frac 1 i  \left( \sum_{j=2}^i  (\Vjk^{\pbet} -\uj^{\pbet}) \uj^{\gamma\beta_*} \frac{\Ejn}{j} \right). 
\]
The first term of the right-hand side is very tedious and delicate to deal with, so we delayed its treatment by stating in Lemma  \ref{Lemmetriplesomme} that it tends to $0$ in probability when $\pbet>1/2$; the proof of this statement is detailed in subsection \ref{demolemmetriplesomme}. Let us then turn to the second term, and prove that it  tends to $0$, and so will $\sqrt{k} I_{2,n}$ as well.   Applying $(\ref{lemmaVjkpartie2})$ with $a=\pbet$, we have, for $\delta >0$ sufficiently small such  that $\epsilon=(\pbet-\delta+\gamma\beta_*)/2$ is positive,  
\[
 \left|\sumideuxk A_{i,n} \frac 1 i  \left( \sum_{j=2}^i  (\Vjk^{\pbet} -\uj^{\pbet}) \uj^{\gamma\beta_*} \frac{\Ejn}{j} \right)\right|  
 \leq 
 O_{\bP}(1) k^{-\epsilon} \  \left(\frac 1 k \sumideuxk |A_{i,n}| u_i^{\delta/2 -1}\right)   \ \left( \frac 1 {k^{3/2 - \epsilon} } \sum_{j=2}^k \uj^{-3/2+\pbet+\gamma\beta_*-\delta} \Ejn \right), 
\]
and we conclude using properties $(\ref{sumdet1})$ and $(\ref{sumEi})$  with $a=1-\delta/2$ as well as  property  $(\ref{sumEibis})$ with $a=3/2-2\epsilon$.
\ztrois

It remains  to consider  the last term $I_{3,n}$ of $\sumideuxk A_{i,n}\Sikbet$, and to prove that it is $o_{\bP}(k_n^{-1/2})$. According to the definition of $C_{j,k,\beta}$ in relation $(\ref{decRFj})$ and using the fact that $\sqrt{k} Y_{n-k,n}^{-\gamma\beta_*}= \sqrt{k}  \left( k/n \right)^{\gamma\beta_*} \left(  Y_{n-k,n} / (n/k) \right)^{-\gamma\beta_*} $ converges (thanks to assumption $(\ref{conditionbiais})$),  we have
\[ 
 \sqrt{k} I_{3,n}  = 
 O_{\bP}(1)  \sumideuxk A_{i,n}  \left(  \frac 1 i\sum_{j=2}^i  \Vjk^{\pbet} (\Vjk^{\gamma\beta_*}-1 )\ksipjsurj  \right) =  O_{\bP}(1) \left(I^{(1)}_{3,n}  - I^{(2)}_{3,n}  + I^{(3)}_{3,n}  - I^{(4)}_{3,n}\right),  
\]
where 
\[ 
 \begar{lll}
  I^{(1)}_{3,n}  &= &   \frac{1}{k+1} \sumideuxk A_{i,n} \left( \frac 1 i  \sum_{j=2}^i  \uj^{\pbet+\gamma\beta_*-1} \ksi'_j   \right)  \zun \\
  I^{(2)}_{3,n}  &= &  \frac{1}{k+1} \sumideuxk A_{i,n} \left(  \frac 1 i \sum_{j=2}^i  \uj^{\pbet-1} \ksi'_j   \right)  \zun \\
   I^{(3)}_{3,n}  &= &   \frac{1}{k+1} \sumideuxk A_{i,n} \left( \frac 1 i   \sum_{j=2}^i  (\Vjk^{\pbet+\gamma\beta_*} -   \uj^{\pbet+\gamma\beta_*}) \uj^{-1} \ksi'_j   \right) \zun \\
    I^{(4)}_{3,n}  &= &  \frac{1}{k+1} \sumideuxk A_{i,n} \left( \frac 1 i   \sum_{j=2}^i  (\Vjk^{\pbet} -   \uj^{\pbet}) \uj^{-1} \ksi'_j   \right) .
\finar \]
Relying on  property $(\ref{lemmaVjkpartie2})$ (stated in Lemma \ref{lemmaVjk}, and applied to $a=\pbet+\gamma\beta_*$) and on the fact that $| A_{i,n} | \leq \Ein +2$, we deduce that, for some given $\delta >0$,  
\[
 |I^{(3)}_{3,n}| \leq    O_{\bP}(1) \ \left(  \frac 1 k \sumideuxk  (\Ein +2)  \ui^{-1+\delta/2}\right)  \ \left( \frac{1}{k^{3/2}}  \sum_{j=2}^k \uj^{ -3/2+\pbet-\delta} \ksi'_j \right) .
\]
Hence, properties $(\ref{sumdet1})$,  $(\ref{sumEi})$ and $(\ref{sumxibis})$  imply that $I^{(3)}_{3,n}$ tends to $0$. Completely similarly, we have $ I^{(4)}_{3,n} =o_{\bP}(1)$. Let us prove that $ I^{(2)}_{3,n}$ also  tends to $0$ ($I^{(1)}_{3,n}$ is handled similarly).  Separating the cases $j <i$ and $j=i$ and using the definition of $ \ksi'_j $ in  relation $(\ref{decXij})$ yield 
\[
 I^{(2)}_{3,n} =  \frac{\gamma}{k+1} \sum_{i=3}^k A_{i,n} \frac 1 i   \sum_{j=2}^{i-1}   \uj^{\pbet-1} \Ejn  +  \frac{\bnk}{k+1} \sum_{i=3}^k A_{i,n} \frac 1 i   \sum_{j=2}^{i-1}   \uj^{\pbet+ \gamma\beta_*-1} \Ejn +  \frac{1}{(k+1)^2}  \sumideuxk A_{i,n}  \ui^{\pbet-2}  \ksi'_i .
\]
The convergence to $0$ of the first (resp. the second)  term  is due to properties $(\ref{EmunE})$ and $(\ref{deltaE})$ with $d=0$ (resp. $d=\gamma \beta_*$) in Lemmas \ref{lemMartingalesEE} and \ref{lemMartingalesEdelta}.   For the third term, we use $| A_{i,n} | \leq \Ein +2$ with Lemma \ref{maxexpo} to write, for some given $\delta >0$, 
\[
 \left| \frac{1}{(k+1)^2}  \sumideuxk A_{i,n}  \ui^{\pbet-2}  \ksi'_j  \right| \leq  O_{\bP}(1) \frac{(\log k)^2}{k^{\delta}} \left(  \frac{1}{k^{2-\delta}} \sumideuxk \ui^{\pbet-2}\right) .
\]
The right-hand side tends to $0$ according to $(\ref{sumdet})$, for $0<\delta<\pbet$. This concludes the proof for the term $\sumideuxk A_{i,n}\Sikbet$. 
\zdeux

\subsubsection{Proof of part $(b)$}

Recall that $B_{i,n}=\frac{1}{\gamma_1} \bnk \ui^{\beta_* \gamma} \Ein$. Since $Z_i=U_H(Y_i)$, using Potter-Bounds $(\ref{PotterBounds})$ for $(\Fbar \circ U_H)U_H^{-\beta} \in RV_{-\pbet}$ and 
working on the event $\calEn $ defined in $(\ref{defEnbeta})$, which satisfies $\lim_{n\tinf}\bP(\calEn)=1$, we have, for $\epsilon >0$ (remind that the sign of $\bnk$ is not known), 
\[
\bnk^{-1} \sumideuxk B_{i,n} S_{i,k}  \leq (1+\epsilon)  \frac{\alpha^{\pbet-\epsilon}}{\gamma_1} \frac{1}{k+1} \sumideuxk   \ui^{\beta_* \gamma} \Ein \left( \frac 1 i  \sum_{j=2}^i \  \uj^{\pbet-1-\epsilon}  \ksi'_j \right).
\]
We are going to prove below that this upper bound,  when multiplied by $\sqrt{k} \bnk$, tends to a quantity arbitrary close to $\bet\lambda$ (for $\epsilon$ small and $\alpha$ close to $1$). A very similar job can be done for the lower bound issued from the application of lower Potter-bounds for $(\Fbar \circ U_H)U_H^{-\beta}$ and from the lower bound in the definition of  $\calEn $, and hence we will have proved that $\sqrt{k} \sumideuxk B_{i,n} \Sikbet$ tends to $\bet\lambda$, as announced. Using  $(\ref{decXij})$ to split $\ksi'_j$  into three parts $\gamma+\gamma(\Ejn-1)+\uj^{\gamma\beta_*}\bnk\Ejn$, we obtain a decomposition of $\sqrt{k} \bnk$ times the upper bound above into three terms $T_{B,n}^{(1)}+T_{B,n}^{(2)}+T_{B,n}^{(3)}$. 
\zdeux

Let us prove that the limit of the first term $ T_{B,n}^{(1)} = \sqrt{k}\bnk (1+\epsilon) \alpha^{\pbet-\epsilon} \frac{p}{k+1} \sumideuxk   \ui^{\beta_* \gamma} \Ein \left( \frac 1 i  \sum_{j=2}^i \  \uj^{\pbet-1-\epsilon}  \right)$, as $n\tinf$, is arbitrarily close to $\bet\lambda$ (taking $\epsilon$ sufficiently small and $\alpha$ sufficiently close to $1$). Indeed, if $\pbet\leq 1$, inequality $(\ref{ineqdik})$ (applied with $a=1-\pbet+\epsilon$) implies that, for some positive constants $C_1$ and $C_2$,  
\[ 
\begar{lll}
\frac{p}{k+1}   \sum_{i=2}^k  \ui^{\gamma \beta_* }  \Ein  \left(\frac 1 i  \sum_{j=2}^i \uj^{\pbet-1-\epsilon}  \right)  & \leq  & \frac{p}{\pbet} \frac{1}{k+1} \sum_{i=2}^k  \ui^{\gamma \beta_* +\pbet-1-\epsilon}  \Ein  - C_1    \frac{p}{(k+1)^{\pbet-\epsilon+1}}  \sum_{i=2}^k  \ui^{\gamma \beta_* -1}  \Ein  \zun \\
 &  \geq  &  \frac{p}{\pbet}  \frac{1}{k+1}  \sum_{i=2}^k  \ui^{\gamma \beta_* +\pbet-1-\epsilon}  \Ein  - C_2   \frac{p}{(k+1)^{\pbet-\epsilon+1}}  \sum_{i=2}^k \ui^{\gamma \beta_* -1}  \Ein  . \
\finar \]
Using $(\ref{sumEi})$ with $a=1-\pbet-\gamma\beta_*+\epsilon$ for the first term and $a=1-\gamma\beta_*$ for the second one, as well as  the fact that  $\bnk$   is equivalent to $-\gamma^2 \beta_* D_* C^{-\gamma\beta_*} \left(  \frac{k+1}{n+1} \right)^{\gamma\beta_*}$,  we obtain via assumption (\ref{conditionbiais}) the desired  limit $\bet\lambda$, by making $\epsilon$ tend to $0$ and $\alpha$ tend to $1$, since $-\gamma^2 \beta_* D_* C^{-\gamma\beta_*}\frac{p}{\pbet} \frac 1{\pbet+\gamma\beta_*}=\bet$.
In the case $\pbet>1$, the treatment is similar, using $(\ref{ineqdik2})$ instead of $(\ref{ineqdik})$ above.  
\zdeux

Secondly, in order to prove that $T_{B,n}^{(2)} =\sqrt{k} \bnk (1+\epsilon) \alpha^{\pbet-\epsilon} \frac{p}{k+1} \sumideuxk   \ui^{\beta_* \gamma} \Ein \left( \frac 1 i  \sum_{j=2}^i \  \uj^{\pbet-1-\epsilon} (\Ejn -1) \right)$ tends to $0$, we   separate the terms $j =i$ (easy to handle and omitted) and $j<i$ : in the latter case,  we use property $(\ref{EmunEbis})$ in Lemma \ref{lemMartingalesEE} (with $\pbet-1-\epsilon$ instead of $\pbet-1$) and the fact that  $\sqrt{k} \bnk$ converges. 
\zdeux

Finally, let us prove that $T_{B,n}^{(3)} = \sqrt{k} \bnk^2 (1+\epsilon) \alpha^{\pbet-\epsilon} \frac{p}{k+1} \sumideuxk   \ui^{\beta_* \gamma} \Ein \left( \frac 1 i  \sum_{j=2}^i \  \uj^{\pbet+\beta_* \gamma-1-\epsilon} \Ein  \right)$ tends to $0$. Using the fact that $\sqrt{k} \bnk^2$  tends to $0$, we bound this term from above by : 
\[
  o_{\bP}(1) \left(  \frac{1}{k+1} \sumideuxk  \ui^{\beta_* \gamma-1} \Ein   \right)  \left(  \frac{1}{k+1} \sumideuxk  \ui^{\pbet+\beta_* \gamma-1-\epsilon} \Ein   \right).
\]
We conclude the treatment  of this term by using $(\ref{sumEi})$.

\subsubsection{Proof of parts $(c)$ and $(d)$}

By the definition of $\Sikbet$ in $(\ref{defSikbeta})$, and the inequality $(\ref{ineqci})$ in Lemma \ref{lem-majorationsdeterministes}, use of Potter-bounds for $(\bar{F} \circ U_H) U_H^{-\beta} \in RV_{-\pbet}$  yields  that, for $\epsilon >0$, 
\[
  \left| \sum_{i=2}^k \deltanmiun  c_i   \Sik \right| 
  \leq  
   (1+\epsilon)  \sum_{i=2}^k \frac{1}{i^2} \  \sum_{j=2}^i  \Vjk^{\pbet-\epsilon} \ksipjsurj . 
\]
Now, working on the event $\calEn $, which satisfies $\lim_{n\tinf}\bP(\calEn)=1$, we have, for $\epsilon >0$ and $\delta >0$, 
\[
 \left|\sum_{i=2}^k \deltanmiun  c_i   \Sikbet \right| 
 \; \leq \;   
 \alpha^{\pbet-\epsilon} (1+\epsilon)  \sum_{i=2}^k \frac{1}{i^2} \ui^{\pbet-\epsilon} \;  \sum_{j=2}^i  \ksipjsurj 
 \; \leq \; 
  cst \left( \frac 1{k^{2-\delta}}  \sum_{i=2}^k \ui^{\pbet-2-\epsilon} \right) \, \left( \frac{1}{k}  \sum_{j=2}^k \uj^{\delta-1} \xi'_j \right) .
\]
Using $(\ref{sumdet})$ and $(\ref{sumxi})$, we see that this expression is lower than $O_{\bP}(1)\times k^{-\pbet+\epsilon+\delta}$, so that part $(c)$ is proved 
as soon as  $\pbet>1/2$, since $\delta$ and $\epsilon$ can be chosen arbitrarily  small. 
\zdeux

Finally, the definition of $\Sikbet$ in $(\ref{defSikbeta})$ on one hand, and  the relation $(\ref{proprieteRnj})$ satisfied by the remainder term $R_{n,i}$ on the other hand, imply that (by inverting sums) 
$$
\left| \sum_{i=2}^k R_{n,i} \Sikbet \right| \leq o_{\bP}(\bnk) \sum_{j=2}^k  \RFjbet \ksipjsurj \log_+ (1/\uj).
$$
As usual, Potter-bounds for $(\bar{F} \circ U_H) U_H^{-\beta} \in RV_{-\pbet}$  yield  that, for $\epsilon >0$, on the event $\calEn$, we have
\[
 \left|\sum_{i=2}^k R_{n,i} \Sik \right| 
 \; \leq \; 
 o_{\bP}(\bnk) \frac{1}{k+1} \sum_{j=2}^k  \uj^{\pbet-1-\epsilon} \ksi'_j  \log_+ (1/\uj). 
\]
Now property $(\ref{sumxi})$ and the fact that $\sqrt{k} \bnk$  converges conclude the proof. 

\subsection{Proof of Proposition \ref{prop-Rn0-and-otherTterms}} \label{subsectionProp2}

\subsubsection{Proof of parts $(a)$ and $(b)$}
\label{partabProp2}

 Concerning the remainder term $R_n^{(0)}$, since  
$\frac{1}{\pbet}=\int_0^1 u^{\pbet-1} du = \sumjdeuxk  \int_{\uj}^{u_{j+1,k}} u^{\pbet-1} du  +  u_{2,k}^{\pbet}/\pbet$,  we obtain
\[ 
\textstyle
 R_n^{(0)}  \; = \; \frac{\gamma}{k+1}\sumjdeuxk \uj^{\pbet-1} - \frac{\gamma}{\pbet} 
\; = \; \gamma    \sumjdeuxk  \int_{\uj}^{u_{j+1,k}} ( \uj^{\pbet-1} -   u^{\pbet-1}) du   -  \frac{\gamma 2^{\pbet}}{\pbet(k+1)^{\pbet}} . 
\]
Using the mean value theorem leads to 
\[
\sqrt{k}  | R_n^{(0)} | \leq   \gamma (1-\pbet) \sqrt{k}   \frac{1}{(k+1)^2} \sumjdeuxk  \uj^{\pbet-2}  + O (k^{1/2-\pbet}).
\]
and we conclude using property $(\ref{sumdet})$ and the condition $\pbet > 1/2$. 
\medskip

 Recall that $R_n^{(1)} =   \sumjdeuxk \RFchapj \frac{\Rnjbet}{j}  $, where  $\Rnjbet$ is defined in $(\ref{defRnjbeta})$. We write $ \RFchapj= \sum_{i=2}^j (\widehat{RF}_{i}- \widehat{RF}_{i-1})$, where we note  $\widehat{RF}_{1} =0$. Hence, inverting sums, we obtain 
\[
R_n^{(1)} =  \sumideuxk (\widehat{RF}_{i}- \widehat{RF}_{i-1})   \sum_{j=i}^k \frac{\Rnj}{j}  \left( \frac{\Znmjun}{\Znmk} \right)^{-\beta} + \frac{\beta}{2}  \sumideuxk  \RFchapj  \frac{\ksi^2_j}{j^2}   \left( \frac{\tilde{Z}_{j,n}}{\Znmk} \right)^{-\beta}, 
\]
where  $\Znmj \leq \tilde{Z}_{j,n} \leq \Znmjun$.

\noindent The definition of $\Fbarn$ implies that $\widehat{RF}_{i}- \widehat{RF}_{i-1} = \widehat{RF}_{i} \frac{\delta_{n-i+2,n}}{i-1}$, for $i > 2$.  Thus, using  $(\ref{proprieteRnj})$ and $\sup_{j\geq 2} \RFchapj/\RFj=O_{\bP}(1)$ (see the proof of part  $(c)$ below for details), we have, if we suppose $\beta \geq 0$ (the case $\beta <0$ is very similar and thus ommited),
\[
|R_n^{(1)}| \leq  O_{\bP}(1) o_{\bP}(\bnk)  \sumideuxk  \frac{RF_{i}}{i-1} \log_+ (\frac{1}{\ui}) + O_{\bP}(1)  \sumjdeuxk \RFj \frac{\ksi^2_j}{j^2}  \left( \frac{Z_{n-j,n}}{\Znmk} \right)^{-\beta}. 
\]
 Now,  using the fact that $Z_i=U_H(Y_i)$,  Potter bounds $(\ref{PotterBounds})$ applied to $\Fbar\circ U_H \in RV_{-p}$ and $U_H^{-\beta} = RV_{-\gamma\beta}$ enable us to write that for any given $\epsilon>0$,   
\[
|R_n^{(1)}| \leq   o_{\bP}(\bnk)  \sumideuxk  \frac{V_{i,k}^{p-\epsilon}}{i-1} \log_+ (\frac{1}{\ui})  +O_{\bP}(1)  \sumjdeuxk \Vjk^{p-\epsilon} \frac{\ksi^2_j}{j^2}  (V_{j+1,k})^{\gamma\beta -\epsilon}.  
\]
Working on the event $\calEn$ which satisfies  $\lim_{n\tinf}\bP(\calEn)=1$, for every $\alpha>1$,  and using the  fact that $\sqrt{k} \bnk$ converges, imply that  
\[
\sqrt{k} |R_n^{(1)}| \leq  o_{\bP}(1)   \frac{1}{k+1} \sumideuxk  \ui^{p-\epsilon-1} \log_+ (\frac{1}{\ui}) + \frac{1}{k^{3/2}}  \sumideuxk  \ui^{\pbet-2\epsilon-2} \ksi^2_j.
\]
We conclude by $(\ref{sumdet1})$ and $(\ref{sumxibis})$ with $\pbet > 1/2$. 

\subsubsection{Proof of part $(c)$}
\label{partcProp2}

Let us now deal with the term $\Tkn{1,2}$, which is defined in relation (\ref{decompTkn1}) and is a delicate part of the proof, and the only one which will require survival analysis arguments. We start by applying the bounds $0\leq -\log(1-x)-x\leq x^2/(1-x)$ $(\forall x<1)$   to $x=1-\RFchapj/\RFj$ for every $j\geq 2$ (which ensures that $\RFchapj>0$ and so $x<1$), yielding
$$
 0 \; \leq \; \Tkn{1,2} \; \leq \; \sumjdeuxk \frac{\RFj}{\RFchapj} \left( 1- \frac{\RFchapj}{\RFj} \right)^2 \RFjbet \ksipjsurj.
$$
We then rely on the so-called Daniels bounds proved in \cite{GillCMS} (page 39) and \cite{Zhou91} (Theorem 2.2), which state that both $\Fbarn(t)/\Fbar(t)$ and its inverse are bounded in probability uniformly for $t<Z_{n,n}$. Since the index $j$ is at least equal to $2$, this implies that $\sup_{j\geq 2} \RFj/\RFchapj=O_{\bP}(1)$. Then (as in the previous subsection \ref{partabProp2}) using the fact that $Z_i=U_H(Y_i)$,  Potter bounds applied to $(\Fbar\circ U_H) U_H^{-\beta} \in RV_{-\pbet}$ enable us to write that for any given $\epsilon>0$,  
$$
 0 \; \leq \; \Tkn{1,2} \; \leq \; O_{\bP}(1) \sumjdeuxk \left( \frac{\RFchapj}{\RFj} - 1 \right)^2  u_{j,k}^{\pbet-1-\epsilon}  (\Vjk/\uj)^{\pbet-\epsilon} \,\ksi'_j .
$$
Now, Theorem 2.1 in \cite{Gill83} applied to the function $h(t)=(\Hbar(t))^{(1+\epsilon)/2}$ guarantees that 
\begeq{theoremGill}
 \sup_{t<Z_{n,n}} \sqrt{n} \, h(t) \left|\frac{\Fbarn(t)-\Fbar(t)}{\Fbar(t)}\right| =O_{\bP}(1),
\fineq
a property which will be applied to $t=Z_{n-j+1,n}$ for every $2\leq j\leq k$ below. Now writing  $\RFchapj/\RFj-1=(\Fbar(\Znmk)/\Fbarn(\Znmk))(W_{n-j+1}-W_{n-k})$ where $W_i=(\Fbarn(Z_{i,n})-\Fbar(Z_{i,n}))/\Fbar(Z_{i,n})$, the combination of the crucial statement (\ref{theoremGill}) with the fact that $h^{-2}$ is nondecreasing, leads to the following bound, working on the set $\calEn$,
$$
 0 \; \leq \; \Tkn{1,2} \; \leq \; O_{\bP}(1) 
 \frac 1 n \frac 1{k+1} \sumjdeuxk  (\Hbar(\Znmjun))^{-1-\epsilon} u_{j,k}^{\pbet-1-\epsilon}  \ksi'_j  . 
$$
Applying then Potter-bounds $(\ref{PotterBounds})$ to the function $(\Hbar^{-1-\epsilon})\circ U_H \in RV_{1+\epsilon}$ then implies that, on the set $\calEn$, we have, for any $\delta>0$,
$$
 0\; \leq \; \sqrt{k} \Tkn{1,2} \; \leq \; O_{\bP}(1)  (\Hbar(\Znmk))^{-1-\epsilon} \left(\frac k n\right) \, k^{-\delta} . \left[ \frac 1{k^{3/2-\delta}} \sumjdeuxk   u_{j,k}^{\pbet-2-3\epsilon}  \ksi_j \right].
$$
First, due to $(\ref{sumxibis})$ in Lemma \ref{lem-weightedsumsbis}, the expression in brackets in the right-hand side of the previous relation is $o_{\bP}(1)$ when $\pbet>1/2$, as soon as $\delta$ and $\epsilon$ are sufficiently small so that  $\pbet>1/2+\delta+ 3\epsilon$. Therefore, since $\Hbar(\Znmk)/(k/n) \stackrel{\bP}{\rightarrow} 1$ as $n\tinf$, all that is left to prove is that $(n/k)^{\epsilon}k^{-\delta}\rightarrow 0$ as $n\tinf$. This is true when assumption (\ref{conditionbiais}) holds with $\lambda\neq 0$, since the latter quantity is equivalent to $\lambda^{-2\delta} (n/k)^{\epsilon-2\gbetaet\delta}$, which indeed converges to $0$ for $\epsilon$ sufficiently small. When assumption (\ref{conditionbiais}) holds with $\lambda=0$, then we use the additional assumption that $n=O(k^b)$ for some $b>1$, which immediately yields $(n/k)^{\epsilon}k^{-\delta}\to 0$ for $\epsilon$ small enough. Part $(b)$ of Proposition \ref{prop-Rn0-and-otherTterms} is thus proved. 

\subsubsection{Proof of part $(d)$}
\label{partdProp2}

Recall that  $\Tkn{1,1,2} =  \sumjdeuxk  \left(   \log\RFj+ \frac{1}{\gamma_1} \log \frac{\Znmjun}{\Znmk}   \right) \RFjbet \ksipjsurj $. Assumption $(\ref{condFbar})$ implies that 
\[
\frac{\Fbar(zt)}{\Fbar(t)} = z^{-1/\gamma_1} \left( 1 + D_1 t^{-\beta_1} (z^{-\beta_1} -1) (1+ o_t (1)) \right) .
\]
Hence, 
\begin{eqnarray*}
 \log\RFj+ \frac{1}{\gamma_1} \log \frac{\Znmjun}{\Znmk}   & =  & \log \left( 1+ D_1 (Z_{n-j+1,n}^{-\beta_1} - Z_{n-k,n}^{-\beta_1}) (1+ o_k (1))  \right)  \zun \\
  & =& D_1  (Z_{n-j+1,n}^{-\beta_1} - Z_{n-k,n}^{-\beta_1}) (1+ o_k (1))   + L_{n,j}, 
\end{eqnarray*}
where $ 0 \leq L_{n,j} \leq D_1^2  (Z_{n-j+1,n}^{-\beta_1} - Z_{n-k,n}^{-\beta_1})^2 (1+ o_k (1))$. Consequently, 
 \begin{eqnarray*}
\Tkn{1,1,2}  & =  &   D_1 (1+\oPdeun) \Znmkbetaun \sumjdeuxk \left(\left(\frac{\Znmjun}{\Znmk}\right)^{-\beta_1} - 1\right) \RFjbet \ksipjsurj 
 \; + \; \sumjdeuxk L_{n,j} \RFjbet \ksipjsurj  \zun \\
  &  = & \Tkn{1,1,2,1}  + \Tkn{1,1,2,2} . 
\end{eqnarray*}

Now, in order to prove that  $\sqrt{k} \Tkn{1,1,2,1}$  tends to $\lambda b_{KM}$, which is defined in the statement of Proposition \ref{prop-Rn0-and-otherTterms}, we deal with the following non-negative quantity, which is equivalent in probability to $\Tkn{1,1,2,1}/(- D_1)$ 
\[
 \widetilde T_{k,n}^{(1,1,2,1)} = \Znmkbetaun \sumjdeuxk \left( 1 - \left(\frac{\Znmjun}{\Znmk}\right)^{-\beta_1} \right) \RFjbet \ksipjsurj.
\] 
Using lower and upper  Potter-bounds for $U^{-\beta_1}_H \in RV_{-\gamma \beta_1}$ and $(\bar{F} \circ U_H) U_H^{-\beta} \in RV_{-\pbet}$ yields, for $\epsilon >0$, 
\begin{eqnarray*}
  & & \sqrt{k}  \Znmkbetaun  \left\{  \frac{1-\epsilon}{k+1} \sumjdeuxk \Vjk^{\pbet-1+\epsilon} \ksi'_j 
  - \frac{(1+\epsilon)^2}{k+1} \sumjdeuxk \Vjk^{\gamma\beta_1 +\pbet-1-2\epsilon} \ksi'_j 
  \right\}  \\ 
  & & \leq \ \sqrt{k}\widetilde T_{k,n}^{(1,1,2,1)}  \ \leq \ 
  \sqrt{k}  \Znmkbetaun   \left\{  \frac{1+\epsilon}{k+1} \sumjdeuxk \Vjk^{\pbet-1-\epsilon} \ksi'_j 
  - \frac{(1-\epsilon)^2}{k+1} \sumjdeuxk \Vjk^{\gamma\beta_1 +\pbet-1+2\epsilon} \ksi'_j 
  \right\} 
\end{eqnarray*}
But $ \Znmkbetaun  = C^{-\gamma\beta_1}(1 + o(1)) (\frac{k}{n})^{\gamma\beta_1}$ (the constant $C$ appears in  formula $(\ref{condUH})$), so $\sqrt{k}  \Znmkbetaun$ tends to $0$ when $\beta_1 > \beta_2$ (due to $(\ref{conditionbiais})$) and,  when $\beta_1 \leq  \beta_2$,  $\sqrt{k}  \Znmkbetaun$ is equivalent  to  $\lambda C^{-\gamma\beta_*} = \lambda C^{-\gamma\beta_1}$.  
Moreover, using $(\ref{defEnbeta})$ with   $\lim_{n\tinf}\bP(\calEn)=1$ and property $(\ref{sumxi})$,  we prove that $\frac 1 k \sumjdeuxk  \Vjk^{\pbet-1+\gamma \beta_1 \pm 2\epsilon} \ksi'_j $ tends to $\frac{\gamma}{\pbet+\gamma \beta_1 \pm 2\epsilon}$ and $\frac 1 k \sumjdeuxk  \Vjk^{\pbet-1 \pm \epsilon} \ksi'_j $ tends to $\frac{\gamma}{\pbet \pm \epsilon}$. After some simplifications, we prove that  $\sqrt{k} \Tkn{1,1,2,1}$  tends to $b_{KM}$, in Porbability,  by making $\epsilon \rightarrow  0$. 
\zun 

Finally, concerning $\Tkn{1,1,2,2} = \sumjdeuxk L_{n,j} \RFjbet \ksipjsurj$, where $0 \leq L_{n,j} \leq D_1^2 ( \Znmjunbetaun - \Znmkbetaun )^2 (1+\oPdeun)$, we use Potter-bounds as previously  to find that , for any given $\epsilon >0$, 
\[
\sqrt{k}  | \Tkn{1,1,2,2} | \leq O(1) \sqrt{k} Z^{-2\beta_1}_{n-k,n}  \ \sumjdeuxk \left( (1+\epsilon)\Vjk^{\gamma \beta_1-\epsilon}  -1 \right)^2 \Vjk^{\pbet-\epsilon}  \ \ksipjsurj
\]
and we proceed as for $ \Tkn{1,1,2,1}$ to prove that $\sqrt{k} \Tkn{1,1,2,2}$ tends to $0$, in Probability.

\subsection{Proof of Proposition \ref{prop-Tkn2Tkn3-mild}} \label{subsectionProp3}

Let us first establish formula $(\ref{lem-decompTkn2Tkn3})$. Recall that (see $(\ref{decompTkn2})$)
\[
\Tkn{2} = \frac{\gamma}{k+1} \sumideuxk \left( \Vjkpbet -\ujpbet \right) \uj^{-1} + \sumideuxk   \Vjkpbet \frac{\ksi'_j -\gamma}{j} +  \sumideuxk \left( \RFjbet - \Vjkpbet \right) \ksipjsurj . 
\]
The definition of $\ksi'_j $ as well  as   decompositions $(\ref{decXij})$ and  $(\ref{decRFj})$ yield 
\begin{align}
\Tkn{2} =  \ &   \  \   \frac{\gamma}{k+1} \sumjdeuxk  \frac{1}{\uj}   (\Vjkpbet-\ujpbet) \Ejn  +     \frac{\gamma}{k+1} \sumjdeuxk  \uj^{\pbet-1} (\Ejn -1)  \nonumber \\
  &  +   \frac{\bnk}{k+1} \sumjdeuxk \uj^{\pbet-1+\gamma\beta_*} \Ejn   +   \frac{\bnk}{k+1} \sumjdeuxk \frac 1{\uj} (\Vjkpbet-\ujpbet) \uj^{\gamma\beta_*}  \Ejn   +  \sum_{j=2}^k \Vjkpbet C_{j,k,\beta} \ksipjsurj  . \nonumber
\end{align} 
The last three terms  of the right-hand side are left unchanged. By applying  decomposition $(\ref{decompVjkpujkp})$ to the first term, we obtain the desired decomposition $(\ref{lem-decompTkn2Tkn3})$. In particular, we can see that the second term of the  right-hand side above vanishes. 
\medskip

 Now, in order to prove the asymptotic result for $\Tkn{2}$,  we  rely of course on the  development $(\ref{lem-decompTkn2Tkn3})$  in 7 different terms.  These terms will be treated separately, one at a time. %
\begitem 
\item[$(a)$]  Concerning the first term, when $\pbet<1$, relation  $(\ref{ineqdik})$   implies that 
\[
 \left| \frac{\pbet  \gamma}{k+1} \sumjdeuxk (\Ejn-1)\left( \frac 1 j \sumideuxj \uipmunbet \!\!- \frac{1}{\pbet} \ujpmunbet \right) \right|  
 \leq O(k^{-\pbet})  \frac 1 k \sumjdeuxk |\Ejn-1| \uj^{-1}
 \leq O(k^{\delta -\pbet}) \, \frac 1 k \sumjdeuxk |\Ejn-1| \uj^{\delta-1}.
 \]
Property $(\ref{sumEi})$ yields that this quantity is  $o_{\bP}(k^{-1/2})$ when $\pbet> 1/2$ for $\delta$ small enough. When $\pbet>1$,  we use $(\ref{ineqdik2})$ instead of $(\ref{ineqdik})$ above.
\smallskip

\item[$(b)$]  Concerning the second term  $\frac{\gamma \pbet}{k+1} \sumjdeuxk  (\Ejn-1)  \ujpmunbet \left( \sum_{i=j}^k \frac{\Ein-1}{i} \right) $,  separating  $i=j$ from $i \geq j+1$ in the sum yields that it is equal to 
\[
 \frac{\gamma \pbet}{(k+1)^2} \sumjdeuxk  (\Ejn-1)^2  \uj^{\pbet-2} 
 + \frac{\gamma \pbet}{k+1} \sum_{i=3}^k (\Ein-1) \left(\frac 1 i \sum_{j=2}^{i-1}   \ujpmunbet  (\Ejn-1) \right). 
\]
Properties $(\ref{sumEibis})$ and $(\ref{EmunEmun}$) prove that this quantity is  $o_{\bP}(k^{-1/2})$ when $\pbet> 1/2$.
\smallskip

\item[$(c)$]  The third term in formula $(\ref{lem-decompTkn2Tkn3})$  is a bias term. Indeed, the expression of $\bnk$ and property $(\ref{sumEi})$ show that 
\[
 \sqrt{k} \bnk \frac 1 {k+1} \sumjdeuxk \uj^{\pbet-1+\gamma\beta_*} \Ejn= -\frac{\gamma^2 \bet D_* C^{-\gamma \bet}}{\pbet + \gamma \bet} (1+o_{\bP}(1))  \sqrt{k} \left( \frac{k}{n} \right)^{\gamma \beta_*} ,
\] 
which yields a  part of the bias term appearing in the statement of Proposition \ref{prop-Tkn2Tkn3-mild}.
\smallskip

\item[$(d)$] The fourth term is   $R_{k,n}= \frac{\bnk}{k+1}  \sumjdeuxk  (\Vjkpbet-\ujpbet)  \uj^{\gamma\beta_*-1}  \Ejn $.
Since $\sqrt{k} b_{n,k} =O(1)$, we have, 
\[
  |R_{k,n}|  \leq  O(1)  \ \sqrt{k} \max_{2\leq j\leq k} \frac{|\Vjkpbet-\ujpbet|}{\uj^{\pbet-1/2-\delta/2}}  \frac{1}{k^2} \sumjdeuxk   \uj^{\pbet+\gamma\beta_*-3/2-\delta/2}  \Ejn ,   
 \]
and  properties  $(\ref{lemmaVjkpartie2})$  and $(\ref{sumEibis})$ imply that $\sqrt{k} R_{k,n}= o_{\bP}(1)$. 
\smallskip

\item[$(e)$] The fifth term  $B_{k,n}=\sum_{j=2}^k \Vjkpbet C_{j,k\beta} \ksipjsurj = (1+\oPdeun) Y_{n-k,n}^{-\gbetaet} D_{\beta} C^{-\gbetaet} \tilde{B}_{k,n}$ will provide the second bias term, where we have noted  
$\tilde{B}_{k,n} = \sum_{j=2}^k \Vjkpbet  (\Vjk^{\gbetaet}-1) \ksipjsurj$,
which is equal to the sum of 2 terms 
\[
 \tilde{B}^{(1)}_{k,n}=   \frac 1 {k+1} \sum_{j=2}^k \uj^{\pbet-1} (\uj^{\gamma \beta_*}-1) \ksi'_j
 \makebox[1.2cm][c]{and} 
 \tilde{B}^{(2)}_{k,n}= \sum_{j=2}^k \left( \Vjkpbet  (\Vjk^{\gbetaet}-1) -  \uj^{\pbet} (\uj^{\gamma \beta_*}-1)\right)\ksipjsurj.
\]
Property $(\ref{sumxi})$ shows that $ \tilde{B}^{(1)}_{k,n}$ converges to $\frac{\gamma}{\pbet+ \gbetaet} - \frac{\gamma}{\pbet} = \frac{-\gamma^2\beta_*}{\pbet(\pbet+\gamma\beta_*)}$. 
On the other hand, we obviously have 
\[
|\tilde{B}^{(2)}_{k,n}|  \leq  \frac 1 {k+1} \sum_{j=2}^k   |\Vjk^{\pbet+\gbetaet} -\uj^{\pbet+\gamma \beta_*}| \uj^{-1} \ \ksi'_j +  \frac 1 {k+1} \sum_{j=2}^k   |\Vjk^{\pbet} -\uj^{\pbet}| \uj^{-1} \ \ksi'_j.  
\]
If we show that $\frac 1 {k+1} \sum_{j=2}^k   |\Vjk^{a} -\uj^{a}| \uj^{-1} \ \ksi'_j$ tends to $0$ for $a=\pbet$ and $a=\pbet+\gbetaet$, then, since $ Y_{n-k,n}^{-\gbetaet}$ is equivalent to $  \left( \frac k n  \right)^{\gbetaet}$,  according to $(\ref{conditionbiais})$ we will have  proved that $B_{k,n} = -  \frac{\gamma^2  D_{\beta} \beta_* C^{-\gbetaet}}{\pbet(\pbet+ \gamma \beta_*)}  \left( \frac{k}{n} \right)^{\gbetaet}  +o_{\bP} ( \frac{1}{\sqrt{k}})$. To do so, we write  
\[
\frac 1 {k+1} \sum_{j=2}^k   |\Vjk^{a} -\uj^{a}| \uj^{-1} \ \ksi'_j  \leq O(1)  \sqrt{k}  \max_{2\leq j\leq k} \frac{|\Vjk^{a}-\uj^{a}|}{\uj^{a-1/2-\delta/2}} \ \frac{1}{k^{3/2}} \sum_{j=2}^k \uj^{a-3/2-\delta/2}  \ksi'_j .
\]
Since $a=\pbet$ or $\pbet+\gbetaet$ are both $> 1/2$, properties $(\ref{sumxibis})$ and $(\ref{lemmaVjkpartie2})$ conclude the proof for the fifth  term. 
\smallskip\\

\item[$(f)$] The absolute value of the sixth term  is shown, thanks to inequality $(\ref{ineq1surj})$, to be lower than
\[
  \pbet \,\gamma  \frac{1}{(k+1)^2}  \sumjdeuxk  \uj^{\pbet-2} |\Ejn|.
\]
Use of  $(\ref{sumEibis})$ with $a=2-\pbet$  and assumption $\pbet> 1/2$ yields that this term is $o_{\bP}(k^{-1/2})$.  
\smallskip

\item[$(g)$] Finally, we deal with the seventh and last term 
$R_{k,n}=\frac{(\pbet)^2\,\gamma}{2(k+1)} \sumjdeuxk \frac 1{\uj} \tilde{V}^{\pbet}_{j,k} \left(  \log(\Vjk/\uj) \right)^2 \Ejn$, where $\tilde V_{j,k}$ lies between $\Vjk$ and $\uj$. On the event $\calEn$, we have
\[
  | R_{k,n}  | 
   \; \leq \; cst \max_{2\leq j\leq k} \frac{|\Ejn|}{\log k}  \ \frac{\log k}{k+1} \sumjdeuxk \frac 1{\uj} \tilde{V}^{\pbet}_{j,k} \left(  \log(\Vjk/\uj) \right)^2 
  \; \leq \;   O_{\bP}(1)  \frac{\log k}{k+1}  \sumjdeuxk\uj^{\pbet-3} (\Vjk-\uj)^2,   
\]
where the mean value theorem and Lemma \ref{maxexpo} were used for the second bound.  
Therefore, for $\delta >0$, 
\[
 | R_{k,n}| \leq o_{\bP}(1)  \  \left( \sqrt{k} \max_{2\leq j\leq k}  \frac{|\Vjk-\uj|}{\uj^{1/2-\delta/2}} \right)^2  \; \frac{k^{\delta}}{(k+1)^2}  \sumjdeuxk\uj^{\pbet-2-\delta}.
\]
and properties $(\ref{sumdet} )$ and $(\ref{lemmaVjkpartie2})$  (with $a=1$) yield $\sqrt{k} R_{k,n} = o_{\bP}(1)$. 
\smallskip
\finit

\subsection{Proof of Lemma \ref{lem-majorationsdeterministes}} \label{subsectionproofLemme1}

Lemma  \ref{lem-majorationsdeterministes} contains a number of different statements, the third and fourth ones being the most relevant in the context of this paper. 
\smallskip

Relation (\ref{ineqci}) is a simple consequence of  the inequality $-x^2\leq \log(1-x)+x \leq 0$ ($\forall x\in [0,1/2]$) applied to $x=1/i$. Then, since $U(j):=\sum_{i=j}^k 1/i = \frac 1{k+1} \sum_{i=j}^k 1/\ui$,  relation (\ref{ineq1surj}) comes from the fact that $\log((k+1)/j)=\sum_{i=j}^k \int_{\ui}^{\uiplusun} x^{-1} dx$ is included in the interval
$$
 \textstyle \left[ \, \frac 1 {k+1} \sum_{i=j}^k 1/{\uiplusun} \; , \; \frac 1 {k+1}\sum_{i=j}^k 1/ {\ui} \, \right] 
 = \left[ U(j) -1/j + \frac 1 {k+1} , U(j) \right]
 \subset \left[ U(j) -1/j , U(j) \right] .
$$

The spirit of the proof of relation (\ref{ineqdik}) is similar : for a given $0<a<1$, setting $\Delta_{i,k}=\ui d_{i,k}$ and noting that $\ui^{1-a}/(1-a)=\int_0^{\ui} u^{-a}du$, we have
\begin{eqnarray*}
\Delta_{i,k} & = & \frac 1{k+1} \sumjdeuxi \uj^{-a} - \frac{1} {1-a} \ui^{1-a} 
 = \sumjdeuxi \int_{\ujmoinsun}^{\uj} (\uj^{-a} -u^{-a})\, du  - u_{1,k}^{1-a}/(1-a)
 \\
 & = & \sumjdeuxi \uj^{1-a} \int_{1-1/j}^1 (1-t^{-a})\,dt - \frac 1{(1-a)(k+1)^{1-a}} 
 \\ 
 & = & \frac 1{(1-a)(k+1)^{1-a}} \left[ \, \sumjdeuxi j^{1-a} \left( \left(1-\frac 1 j\right)^{1-a} - \left(1-\frac {1-a}{j}  \right)  \right)  \; - \; 1 \, \right]  
\end{eqnarray*}
Applying, for each $j$, the Taylor formula of order 2 to the function $x\to (1-x)^{1-a}-(1-(1-a)x)$ between $0$ and $1/j$ (which is lower than $1/2$) leads to the following bounds 
$$
 -1 - a(1-a) 2^{a} \sumjdeuxi j^{-1-a} 
\; \leq \; (1-a)(k+1)^{1-a} \Delta_{i,k} \; \leq \; 
 -1 - \frac {a(1-a)}{2} \sumjdeuxi j^{-1-a} 
$$
and therefore we have shown that, when $0<a<1$,  statement (\ref{ineqdik}) holds for instance with the constants $C_1=1/(1-a)$ and $C_2=[1+a(1-a)2^a(\zeta(1+a)-1)]/(1-a)$.  This means in particular that the values $d_{i,k}$ are always negative, which is a  fact often used in the proofs of this paper.
\smallskip

The proof of (\ref{ineqdik2}) when $a<0$ is performed similarly : we come up to 
$$
 d_{i,k} = - \frac 1{\ui (1-a)(k+1)^{1-a}} - \frac {a}{2(k+1)^{-a}} \frac 1 i \sumjdeuxi j^{-1-a} (1-c_j)^{-1-a}
$$
where $c_j$ are values between $0$ and $1/j$ for each $2\leq j\leq i$ (thus lower than $1/2$). The second term in the right-hand side of the formula above being positive, and since $(k+1)^{1-a}>k+1$, we have proved the lower bound for $d_{i,k}$. For the upper bound, we bound the right-hand side above by zero plus the positive value $(-a/(k+1)) \frac 1 i \sumjdeuxi \uj^{-1-a}$. Distinguishing the cases $a<-1$, $a=-1$ and $-1<a<0$ then leads easily to the desired upper bound.

\subsection{Proof of Lemma \ref{lemMartingalesEE}}   

We first deal with (\ref{EmunEmun}). Letting  $W_{j-1}$ denote $\frac 1 j \sum_{i=2}^{j-1} \uipmunbet (E_i-1)$, we remark that $E_j-1$ and $W_{j-1}$ are  independent and centered, and it is easy to check that the products $(E_j-1)W_{j-1}$ ($j=3\ldots k$) are then centered and uncorrelated. Therefore, it suffices to prove that $\frac 1 {k} \sum_{j=3}^k \bE(W_{j-1}^2)$  (which is equal to the variance of the left-hand side of (\ref{EmunEmun})) converges to $0$. By construction, $\bE(W_{j-1}^2)=\frac 1 {j^2} \sum_{i=2}^{j-1} u_{i,k}^{2(\pbet-1)}$. If $\pbet \leq 1$,  by using the inequality (\ref{ineqdik}) with $a=2(1-\pbet) \in [0,1[$, we have  $\frac 1 {j} \sum_{i=2}^{j-1} u_{i,k}^{2(\pbet-1)}\leq  \frac 1 {2\pbet-1} u_{j-1,k}^{2(\pbet-1)}$.  If $\pbet >1$, we have simply (via $u_{i,k} \leq u_{j-1,k}$) $\frac 1 {j} \sum_{i=2}^{j-1} u_{i,k}^{2(\pbet-1)}\leq   u_{j-1,k}^{2(\pbet-1)}$. We can thus deduce that $\bE(W_{j-1}^2) \leq \frac{cst}{j}  u_{j-1,k}^{2(\pbet-1)} \leq \frac {cst} {k} u_{j-1,k}^{2\pbet-3}$. Finally, we obtain that our quantity of interest $\frac 1 {k} \sum_{j=3}^k \bE(W_{j-1}^2)$ is  lower than a constant times $\frac 1 {k^2} \sum_{j=3}^{k} u_{j-1,k}^{2\pbet-3}$, which converges to $0$ because $\pbet>1/2$.
\medskip

Concerning (\ref{EmunE}), defining now $W_{j-1}$ as $\frac 1 j \sum_{i=2}^{j-1} u^{\pbet+d-1}_i E_i$, the difference with the previous case is that $W_{j-1}$ is not centred. However the products $(E_j-1)W_{j-1}$ are still uncorrelated, and it again suffices to prove the convergence to $0$ of the variance of the left-hand side of (\ref{EmunE}), which is now equal to $\frac 1 {k^2} \sum_{j=3}^k \bE(W_{j-1}^2)$.  By the Cauchy-Schwarz inequality, we have here 
$$
\bE(W_{j-1}^2) \leq \bE\left[ \left( \frac 1 j \sum_{i=2}^{j-1} u_{i,k}^{2(\pbet-1)} \right) 
 \left( \frac 1 j \sum_{i=2}^{j-1} E_i^2 \right) \right] \leq \frac 2 {j} \sum_{i=2}^{j-1} u_{i,k}^{2(\pbet-1)} 
 \leq cst \ u_{j-1,k}^{2(\pbet-1)}
$$
where the last inequality was shown in the treatment of (\ref{EmunEmun}) above.  Therefore, we deduce that $\frac 1 {k^2} \sum_{j=3}^k \bE(W_{j-1}^2)$ is lower than a constant times 
$\frac 1 {k^2} \sum_{j=2}^{k-1} u_{j,k}^{2(\pbet-1)}$, which is $O(k^{-1})$ since $\pbet>1/2$.
\medskip 

Concerning (\ref{EmunEbis}), we invert the two sums and then, we have to deal with 
\[
\frac{1}{k^2} \sum_{j=2}^k  \uj^{\pbet-1} (E_j-1)  \ \left\{  \sum_{i=j+1}^{k} \ui^{\gamma\beta_* -1} E_i \right\} .
\]
Defining now $W_{j+1}$ as $ \sum_{i=j+1}^{k} \ui^{\gamma\beta_* -1} E_i $ which is independent of $E_j-1$, it is easy to check that $(E_j-1)W_{j+1}$ ($j=2\ldots k$) are then centred and uncorrelated. Therefore, it suffices to prove the convergence to $0$ of the variance of the left-hand side of $(\ref{EmunEbis})$, which is equal to $\frac{1}{k^4}   \sum_{j=2}^k \uj^{2(\pbet-1)} \bE(W_{j+1}^2) $. By the Cauchy-Schwarz inequality, we have
\[
 \bE(W_{j+1}^2) \leq  2 (k-j) \sum_{i=j+1}^{k}  \ui^{2(\gamma\beta_* -1)}  \leq  k\sum_{i=j+1}^{k}  \ui^{2(\gamma\beta_* -1)} .
\]
Inverting the two sums we deduce that $\frac{1}{k^4}   \sum_{j=2}^k \uj^{2(\pbet-1)} \bE(W_{j+1}^2) $ is lower than $\frac{1}{k^3}   \sum_{i=3}^k  \ui^{2 \gamma\beta_* -2}  \left( \sum_{j=2}^{i-1}  \uj^{2(\pbet-1)} \right) $, which is lower than $\frac{cst}{k^3}   \sum_{i=3}^k  \ui^{2 \gamma\beta_* -2}  i   \ u_{i,k}^{2(\pbet-1)} \leq \frac{cst}{k^4}   \sum_{i=3}^k  \ui^{-4+\epsilon}$ ($\epsilon >0$), which converges to $0$. 
\zdeux

Concerning finally (\ref{EmunEmunE}), the method developed above works similarly. By noting $W'_{l,n}=l^{-1} \sum_{j=2}^{l-1} u_{j,k}^{\pbet-1} E_j$ and $W_{i,n} = i^{-1} \sum_{l=3}^{i-1} (E_l-1) W'_{l,n}$, the variables $W'_{l,n}$ are not centred but their variance can be shown to be lower than a constant times $u_{l,k}^{2(\pbet-1)}$ . Since $W'_{l,n}$ and $E_l-1$ are independent, the variables $(E_l-1)W'_{l,n}$ are centred and uncorrelated, and  thus $W_{i,n}$ has a variance lower than a constant times $k^{-1} \ui^{-1-2(1-\pbet)}$, and is independent of $(E_i-1)$, so the variance of the left-hand side of $(\ref{EmunEmunE})$ is lower than a constant times $k^{-2} \sum_{i=4}^k \ui^{-1-2(1-\pbet)}$, where $1+2(1-\pbet) < 2$ when $\pbet>1/2$. The proof is then over via relation $(\ref{sumdet})$. 

\subsection{Proof of Lemma \ref{lemMartingalesEdelta} }

First, let us recall that $\delta_i=\bI_{U_i\leq p(Z_i)}$ with $(U_i)$ uniformly distributed and independent of the $Z_i$'s. Then, let us settle the following notations. First, the difference $\deltanmiun -p$ will be systematically cut in three terms
$$
\deltanmiun -p \stackrel{d}{=} \Delta_i^{(1)} +  \Delta_i^{(1)} +  \Delta_i^{(1)} 
\makebox[1.5cm][c]{where} \left\{ 
\begar[c]{rcl}  \Delta_i^{(1)} & = & \bI_{U_i  \leq p}  -p   ,  \\
\Delta_i^{(2)} & = &  \bI_{U_i  \leq p\circ U_H(n/i)} -  \bI_{U_i  \leq p}, \\
\Delta_i^{(3)} & = &   \bI_{U_i \leq p\circ U_H(Y_{n-i+1,n})} -  \bI_{U_i \leq p\circ U_H(n/i)}.
\finar \right.
$$ 
The first of these terms will be the less negligible  one, but the easiest to deal with. The second one will still be simple to handle, but leads to non-centered factors. The third one, $\Delta_i^{(3)}$, will be the "smallest", but the most difficult to deal with, since it is correlated with the observations $(Z_i)$ (and therefore with the variables $\Ejn$). In the sequel, $cst$ will design an absolute positive constant which varies from line to line. 
\zdeux

We start by proving $(\ref{Edelta})$. Setting  $W_{in} = \frac 1 i \sum_{j=2}^{i-1} \ujpmunbet (\Ejn-1) $ and $A^{(m)}_n = k^{-1/2}\sum_{i=3}^k \Delta_i^{(m)} W_{in}$, we intend to prove that $\bV(A^{(1)}_n)$ and $\bV(A^{(2)}_n)$ go to $0$ as $n\tinf$, and that $A^{(3)}_n$ converges to $0$ in probability. Concerning first $A^{(1)}_n$, we note that the variables $\Delta_i^{(1)}$ are i.i.d. centered and independent of the variables $(\Ejn)$ and thus of the centered $W_{in}$ : therefore, the product $\Delta_i^{(1)}W_{in}$ is centered and uncorrelated with $\Delta_{i'}^{(1)}W_{i'n}$ for any $i\neq i'$, and consequently 
$$
{\bV}(A_n^{(1)}) = \frac 1 k \sum_{i=3}^k {\bV}(\Delta_i^{(1)}W_{in}) = \frac 1 k \sum_{i=3}^k p(1-p) {\bV}(W_{in}) \leq \frac {cst}{k^2} \sum_{i=3}^k \ui^{-1-2(1-\pbet)} \tqdninf 0
$$
because $1+2(1-\pbet)<2$ when $\pbet>1/2$.   Above we have bounded ${\bV}(W_{in})$ with similar tools as those used in the proof of Lemma \ref{lemMartingalesEE}, by a constant times 
 $\frac{1}{k+1} \ui^{-1-2(1-\pbet)}$.
 
\smallskip

Concerning now $A_n^{(2)}$, we note that the variables $\Delta_i^{(2)}$ are not centered but are still independent, and independent of the $W_{in}$. Since $W_{in}$ is centered, the products $\Delta_i^{(2)}W_{in}$ are still centered but are now correlated, since, for $i'<i$,
$$
Cov(\Delta_i^{(2)}W_{in},\Delta_{i'}^{(2)}W_{i'n}) 
= \bE(\Delta_i^{(2)})\bE(\Delta_{i'}^{(2)})\bE(W_{in}W_{i'n}) \neq 0.
$$
Using relation $(\ref{pHall})$ of Lemma \ref{fonctionp}, both the variance and the absolute value of the expectation of $\Delta_i^{(2)}$ turn out to be lower than $cst\, (i/n)^{\gamma\beta_*}$, which, due to assumption $(\ref{conditionbiais})$,    
is itself lower than $cst\, k^{-1/2}$. On the other hand, we have, for $i'<i$, 
$$
Cov(W_{in},W_{i'n}) 
= \bE(W_{i'n}W_{in}) 
= \frac {i'} {i} \bE\left( W_{i'n}^2 \right) + \bE\left( W_{i'n} . \frac 1 i \sum_{j=i'}^{i-1} \uj^{\pbet-1}(\Ejn-1) \right) = \frac {i'} {i} \bE\left( W_{i'n}^2 \right) 
$$
Therefore, we may write that (using the bound $\bE(\Delta_i^{(2)})\leq cst\, k^{-1/2}$ in the second term below, but simply bounding $|\Delta_i^{(2)}|$ by $1$ in the first term)
\begin{eqnarray*}
{\bV}(A_n^{(2)}) 
& = & \frac 1 k \sum_{i=3}^k \bE\left((\Delta_i^{(2)})^2 W_{in}^2\right) 
  + \frac 2 k \sum_{i=4}^k \sum_{i'=3}^{i-1} \bE(\Delta_i^{(2)})\bE(\Delta_{i'}^{(2)}) Cov(W_{in},W_{i'n}) \\
& \leq &  \frac 1 k \sum_{i=3}^k \bE\left( W_{in}^2\right) 
  + \frac {cst} {k^2} \sum_{i=4}^k \sum_{i'=3}^{i-1} \bE\left( W_{i'n}^2\right)
\ \leq \  \frac {cst} {k^2} \sum_{i=3}^k \ui^{-1-2(1-\pbet)}
\end{eqnarray*}
and this converges to $0$ when $\pbet>1/2$, as desired. 
\smallskip

In order to finish the proof of $(\ref{Edelta})$, we have to justify that the last part, $A_n^{(3)}$, converges to $0$ in probability. Our proof is based on the important fact that, for any value $\tilde p\in ]1/2,\pbet]$, 
\begeq{relationclef-ecartdeltani}
\frac 1{\sqrt{k}} \sum_{i=3}^{k} | \Delta_i^{(3)} | \; \ui^{\tilde p -1} = 
\frac 1{\sqrt{k}} \sum_{i=3}^{k} \left| \bI_{U_i \leq p\circ U_H(Y_{n-i+1,n})} -  \bI_{U_i \leq p\circ U_H(n/i)} \right| \; \ui^{\tilde p -1} 
\stackrel{\bP}{\longrightarrow} 0.  
\fineq
This result is very close to the one stating that $\sqrt{k}B^{(1)}_{k,n}=o_{\bP}(1)$ in subsection \ref{preuve-morceauprincipal-sommeAinSik}, it is proved completely similarly, therefore details are omitted.    Therefore, in view of relation $(\ref{relationclef-ecartdeltani})$, convergence in probability to $0$ of $A_n^{(3)}$ will follow from the following statement : for every $A>0$, 
\begeq{inegalitemaximaleWin}
\bP \left( \max_{3\leq i\leq k} | W_{in} / \ui^{\tilde p-1} | > A \right) \tqdninf 0.
\fineq 
Considering the sum of independent variables $S_i=\sum_{j=2}^{i-1} j^{\pbet-1}(E_j-1)$ (where $E_j$ denote iid standard exponential variables), we have $W_{in}/\ui^{\tilde p-1} \stackrel{d}{=} k^{\tilde p - \pbet} S_i/i^{\tilde p}$, and therefore, application of the H\'ajek-R\'enyi maximal inequality (see for instance Section 7.4 of \citet{ChowTeicher1997}) leads to 
$$
\bP \left( \max_{3\leq i\leq k} | W_{in} / \ui^{\tilde p-1} | > A \right)
\leq (Ak^{\pbet-\tilde p})^{-2} \sum_{i=2}^k \frac{ \bE( (i^{\pbet-1} (E_i-1) )^2 ) } {i^{2\tilde p}} 
\; = \; \frac 1 {A^2} k^{-2(\pbet-\tilde p)} \sum_{i=2}^k i^{-2+2(\pbet-\tilde p)} 
$$
which goes to $0$ as $n\tinf$, since $0<\pbet-\tilde p<1/2$, and this proves (\ref{inegalitemaximaleWin}). This ends the justification of relation (\ref{Edelta}). 
\medskip

Concerning now relation $(\ref{deltaE})$, we again divide $\delta_{n-i+1,n}-p$ in three parts as above, and the $\Delta_i^{(3)}$ part is proved by combining relation (\ref{relationclef-ecartdeltani}) with Lemma \ref{maxexpo}; the other two parts are easy to deal with. 
\medskip

Concerning relation $(\ref{deltaEmunE})$, we proceed similarly as for $(\ref{Edelta})$, defining now 
$$ 
W_{in}=\frac 1 i \sum_{j=3}^{i-1} (\Ejn-1) \left( \frac 1 j \sum_{l=2}^{j-1} \uj^{\pbet-1} \Eln \right)
\makebox[1.2cm][c]{ and } 
A_n^{(m)}=k^{-1/2}\sum_{i=4}^k \Delta_i^{(m)}W_{in}
\hspace{0.3cm}\mbox{ for $m=1,2,3$.}
$$ 
These variables  $W_{in}$ are still centered, and their variance and covariances can be bounded in exactly the same way as were those of $\frac 1 i \sum_{j=2}^{i-1} (\Ejn-1)\uj^{\pbet-1}$ : therefore, convergence to $0$ of the variances of the corresponding terms $A_n^{(1)}$ and $A_n^{(2)}$ is proved as above. And since $W_{in}$ also possesses an appropriate martingale structure to which the H\'ajek-R\'enyi maximal inequality can be applied, convergence in probability to $0$ of $A_n^{(3)}$ holds, and so does $(\ref{deltaEmunE})$. 
\medskip

Concerning finally relation $(\ref{EmunEdeltamp})$, we write its left-hand side as the sum of the following three expressions, noting $W'_{ln} = \frac 1 l \sum_{j=2}^{l-1} \uj^{\pbet-1}\Ejn$,
$$
A_n^{(1)} \; = \; 
\frac 1{\sqrt{k}} \sum_{i=4}^k (\Ein-1) \left\{ \frac 1 i \sum_{l=3}^{i-1} \Delta_i^{(1)} W'_{ln}  \right\} 
\makebox[1.5cm][c]{,}
A_n^{(2)} \; = \;
\frac 1{\sqrt{k}} \sum_{i=4}^k (\Ein-1) \left\{ \frac 1 i \sum_{l=3}^{i-1} \Delta_i^{(2)} W'_{ln}  \right\} 
$$
and 
$$
A_n^{(3)} \; = \;
\frac 1{\sqrt{k}} \sum_{i=3}^{k-1} \Delta_i^{(3)} \left( \frac 1 i \sum_{j=2}^{i-1} \uj^{\pbet-1}\Ejn \right) \left( \sum_{l=i+1}^k \frac{\Eln-1}{l} \right).
$$
As sums of centered and uncorrelated terms, the quantities $A_n^{(1)}$ and $A_n^{(2)}$ can be handled similarly as previously (with a bit more efforts for $A_n^{(2)}$), and their variances shown to go to zero. Concerning now $A_n^{(3)}$, setting $\tilde S_i= \sum_{l=i+1}^k (\Eln-1)/l$ \ and \ $W_{in}=\frac 1 i \sum_{j=2}^{i-1} \ujpmunbet (\Ejn-1)$, we have, for $\tilde p\in ]1/2,\pbet[$, 
$$
|A_n^{(3)}| \; \leq \; \max_{3\leq i \leq k-1} \ui^{1-\tilde p} |W_{in} \tilde S_i| . \frac 1{\sqrt{k}} \sum_{i=3}^{k-1} |\Delta_i^{(3)}| \ui^{\tilde p-1} 
\; + \;
\max_{3\leq i \leq k-1}  |\tilde S_i| . \frac {cst} {\sqrt{k}} \sum_{i=3}^{k-1} |\Delta_i^{(3)}| \ui^{\pbet-1}  
$$
In view of statements (\ref{relationclef-ecartdeltani}) and (\ref{inegalitemaximaleWin}), we thus have to prove that $\max_{i\leq k} |\tilde S_i|$ is bounded in probability . But since $\max_{i\leq k} |\tilde S_i|\leq  |S_k| + \max_{i\leq k} |S_i|$ where $S_i=\sum_{l=1}^i (\Eln-1)/l$, and $\bV(S_k)=\sum_{l=1}^k 1/l^2\leq \pi^2/6$, the Markov inequality and the usual maximal inequality of Kolmogorov yield the desired result, for any $A>0$, $ \bP[ \max_{3\leq i \leq k-1}  |\tilde S_i| >A ] \leq 8\bV(S_k)/A^2 \leq cst/A^2$, which is as small as desired.

\subsection{Proof of Lemma \ref{Lemmetriplesomme}}  \label{demolemmetriplesomme}

Formula  $(\ref{decompVjkpujkp})$  yields 
\[
\sumideuxk A_{i,n} \frac 1 i  \left( \sum_{j=2}^i  (\Vjkpbet -\ujpbet) \frac{\Ejn}{j} \right) = R_{1,n} + R_{2,n} + R_{3,n}, 
\]
with 
\[ \begar{lll}
R_{1,n} &= &    -\pbet  \sumideuxk A_{i,n} \frac 1 i  \sum_{j=2}^i  \ujpbet \left( \summ{l=j}{k}  \frac{(\Eln-1)}{l} \right)  \frac{\Ejn}{j}  \zun \\
 R_{2,n} &= & -\pbet \sumideuxk A_{i,n} \frac 1 i  \sum_{j=2}^i  \ujpbet \left( \summ{l=j}{k} \frac 1 l -  \log \frac{k+1}{j} \right)  \frac{\Ejn}{j} \zun \\
  R_{3,n} &= &  \frac{\pbet^2}{2}  \sumideuxk A_{i,n} \frac 1 i  \sum_{j=2}^i  \tilde{V}_{j,k}^{\pbet} \left( \log \frac{\Vjk}{\uj} \right)^2 \frac{\Ejn}{j} , 
\finar \]
where $\tilde V_{j,k}$ lies between $\Vjk$ and $\uj$.  The main term is $R_{1,n}$, but we consider $R_{2,n}$ and $R_{3,n}$ first. Inequality $(\ref{ineq1surj})$ in Lemma \ref{lem-majorationsdeterministes} implies that, for $\delta >0$,  
\[
\sqrt{k} |R_{2,n}| \leq O(1) \left( \frac 1 k  \sumideuxk |A_{i,n}| \ui^{\delta -1}\right)   \ \left(\frac{1}{k^{3/2}}   \sum_{j=2}^k \uj^{\pbet-2-\delta} \Ejn \right). 
\]
Hence $ \sqrt{k} R_{2,n}$ tends to $0$ thanks to $(\ref{sumEi})$ and $(\ref{sumEibis})$, with $\pbet> 1/2$. Now, concerning $R_{3,n}$, we proceed as in the proof of Proposition \ref{prop-Tkn2Tkn3-mild} part $(g)$. Using the mean value theorem, Lemma  \ref{maxexpo}  and then applying property $(\ref{lemmaVjkpartie2})$ (with $a=1$),  then, working on  the event $\calEn$ defined in $(\ref{defEnbeta})$, we have, for $\delta >0$, 
\[
\sqrt{k} |R_{3,n}| \leq  o_{\bP} (1)  \left( \frac 1 k  \sumideuxk |A_{i,n}| \ui^{\delta -1}\right)   \ \left(\frac{k^{\delta}}{k^{3/2}}   \sum_{j=2}^k \uj^{\pbet-2-2\delta}  \right),  
\]
and we conclude using $(\ref{sumEi})$ and $(\ref{sumdet})$. We thus have to deal with the first term $R_{1,n}$, and we start by separating the cases  $l=j$ and $l>j$ to obtain  
\[
R_{1,n}  = -\pbet  \sumideuxk A_{i,n} \frac 1 i  \sum_{j=2}^i  \ujpbet \frac{(\Ejn-1)}{j}   \frac{\Ejn}{j}
-\pbet  \sum_{i=2}^{k} A_{i,n} \frac 1 i  \sum_{j=2}^i  \ujpbet \left( \summ{l=j+1}{k}  \frac{(\Eln-1)}{l} \right)  \frac{\Ejn}{j}.
\]
We prove easily that the first term of the right-hand side is $o_{\bP}(1/\sqrt{k})$, using $(\ref{sumEi})$ and $(\ref{sumEibis})$. For the second term, we separate the cases $j=i$ and $j<i$ and obtain 
\[
\sum_{i=2}^{k-1} A_{i,n}  \frac 1 i    \ui^{\pbet}  \left( \summ{l=i+1}{k}  \frac{(\Eln-1)}{l} \right)  \frac{\Ein}{i}
+   \sum_{i=3}^{k} A_{i,n}   \frac 1 i    \sum_{j=2}^{i-1}  \ujpbet  \frac{\Ejn}{j} \left( \summ{l=j+1}{k}  \frac{(\Eln-1)}{l} \right) .
\]
We prove easily, using $(\ref{sumEi})$ and $(\ref{sumEibis})$,  that the first term of the right-hand side is $o_{\bP}(1/\sqrt{k})$. The second term is split in two by separating  the cases $j+1 \leq l \leq i$ and  $i+1 \leq l \leq k$.  We obtain the following two terms 
\[ \begar{lll} 
R'_{1,n} & =  & \sum_{i=3}^{k} A_{i,n}   \frac 1 i    \sum_{j=2}^{i-1}  \ujpbet  \frac{\Ejn}{j} \left( \summ{l=j+1}{i}  \frac{\Eln-1}{l} \right)  \zun \\
 R'_{2,n}  & =&   \sum_{i=3}^{k-1} A_{i,n}   \frac 1 i    \sum_{j=2}^{i-1}  \ujpbet  \frac{\Ejn}{j} \left( \summ{l=i+1}{k}  \frac{\Eln-1}{l} \right).
\finar \]
Inverting the sum in $i$ and the sum in $l$,  we see that  $\sqrt{k} R'_{2,n} $ tends to $0$ thanks to properties $(\ref{EmunEmunE})$ and $(\ref{EmunEdeltamp})$ in Lemmas  \ref{lemMartingalesEE} and \ref{lemMartingalesEdelta}. Now, inverting the sum in $l$ and the sum in $j$ yields 
\[
R'_{1,n} = \sum_{i=3}^{k} A_{i,n}   \frac 1 i    \sum_{l=3}^{i} \frac{\Eln-1}{l}   \left(     \sum_{j=2}^{l-1}  \ujpbet \frac{\Ejn}{j}  \right). 
\]
Separating finally the cases $l=i$ and $l<i$, we obtain the  following two  terms : 
\[ \begar{lll} 
R''_{1,n} & =& \sum_{i=3}^{k} A_{i,n}    \frac{\Ein-1}{i^2}   \left(     \sum_{j=2}^{i-1}  \ujpbet \frac{\Ejn}{j}  \right),    \zun \\
R''_{2,n} & =&\sum_{i=4}^{k} A_{i,n}   \frac 1 i    \sum_{l=3}^{i-1} \frac{\Eln-1}{l}   \left(     \sum_{j=2}^{l-1}  \ujpbet \frac{\Ejn}{j}  \right).
\finar \]
$\sqrt{k}  R''_{2,n}$ tends to $0$ thanks to properties $(\ref{EmunEmunE})$ and  $(\ref{deltaEmunE})$ in Lemmas \ref{lemMartingalesEE} and \ref{lemMartingalesEdelta}. We now  conclude the  proof of this lemma by proving that  $ \sqrt{k}  R''_{1,n}$ tends to $0$.  Since $|A_{i,n} | \leq \Ein +2$, 
\[
\bE(|R''_{1,n}|) \leq  O(1) \frac 1 k  \sum_{i=3}^{k}  \frac{1}{i^2}   \sum_{j=2}^{i-1}  \uj^{\pbet-1},
\]
and the right-hand side tends to $0$ using  $(\ref{ineqdik})$ (or $(\ref{ineqdik2})$) and $(\ref{sumdet})$.

\subsection{Elements of proof for the other lemmas} 
\label{subsectionProofLemmesMineurs}

Concerning Lemmas \ref{lem-weightedsums} and \ref{lem-weightedsumsbis}, relation $(\ref{sumdet1})$ is just the convergence of a Riemann sum, (\ref{sumdet}) is just one definition of the Zeta function, statements (\ref{sumxi}) and (\ref{sumxibis}) have been proved in Lemma  2 of \citet{WormsWorms14} respectively for $0<a<1$ and $a>1$ (for (\ref{sumxi}), the treatment of the case $a\leq 0$ is similar).  Property (\ref{sumEi}) is a  simple application of the triangular law of large numbers, whereas property (\ref{sumEibis})  is deduced easily from  $(\ref{sumdet1})$.  Details are omitted. 
\smallskip

Lemma \ref{maxexpo} is a simple consequence of the fact that the exponential distribution admit a finite exponential moment. Proof of  Lemma \ref{fonctionp} is omitted (see \citet{BeirlantBardoutsosWetGijbels16} for (\ref{pHall})). 
\smallskip

 Lemma \ref{lemmaVjk} is based on the fact that the  uniform empirical quantile process based on a uniform sample of size $k$ satisfies  $\sqrt{k} \sup_{1/(k+1) \leq t \leq k/(k+1)} \left| (\Gamma_k^{-1}(t) -t) / t^{1/2-\delta/2} \right| = O_{\bP} (1)$ (see, for example, \cite{ShorackWellner86} sections 10.3 and 11.5).  Since $\Gamma_k^{-1}(t)= \Vjk$, for $\frac{j-1}{k}  \leq t \leq  \frac{j}{k} $, this yields relation $(\ref{lemmaVjkpartie2})$ for  $a=1$. From the mean value theorem and working on the event $ \calEn$ defined in $(\ref{defEnbeta})$, relation $(\ref{lemmaVjkpartie2})$ for a general $a>0$ follows easily.

\end{document}